\documentclass[11pt,oneside]{amsart}
\usepackage{amsmath}
\usepackage{amsfonts}
\usepackage{amssymb}
\usepackage{amsthm}
\usepackage[colorlinks]{hyperref} 
\usepackage{wasysym}
\usepackage{graphicx}
\usepackage{dashrule}
\usepackage[all]{xy}
\usepackage{multirow}
\usepackage{enumitem}
\usepackage{tikz}
\usetikzlibrary{arrows,decorations.pathmorphing,backgrounds,positioning,fit,calc}
\usepackage{tikz-cd}
\usepackage{framed}

\usepackage{latexsym}
\usepackage{stix}

\normalfont\upshape

\theoremstyle{plain}
\newtheorem{defn}{Definition}[section]
\newtheorem{thm}[defn]{Theorem}
\newtheorem{propn}[defn]{Proposition}

\numberwithin{equation}{section}

\theoremstyle{remark}
\newtheorem{rmk}[defn]{Remark}

\DeclareMathOperator{\GL}{GL}

\DeclareMathOperator{\spec}{Spec}

\DeclareMathOperator{\conv}{Conv}
\DeclareMathOperator{\stab}{Stab}

\newcommand{\nocontentsline}[3]{}
\newcommand{\tocless}[2]{\bgroup\let\addcontentsline=\nocontentsline#1{#2}\egroup}

\newcommand{\CC}{\mathbb{C}}
\newcommand{\PP}{\mathbb{P}}

\newcommand{\GG}{\mathbb{G}}

\newcommand{\dblslash}{/\! \!/}

\newcommand{\env}{/\!/ } 

\newcommand{\inenv}{\dblslash \!_{\circ}}

\newcommand{\ten}{\otimes}

\newcommand{\mf}{\mathfrak}

\newcommand{\kk}{\Bbbk}

\newcommand{\ssfg}{\mathrm{ss,fg}}

\newcommand{\rms}{\mathrm{s}}

\newcommand{\ximp}{X_{{\rm impl}}}
\newcommand{\liets}{{\mathfrak t}^*}
\newcommand{\Proj}{{\rm Proj}}
\newcommand{\Spec}{{\rm Spec}}
\newcommand{\nc}{\newcommand}
\nc{\bla}{\phantom{bbbbb}}

\newcommand{\beq}{\begin{equation}}
\newcommand{\eeq}{\end{equation}}
\newcommand{\barr}{\begin{array}}
\newcommand{\earr}{\end{array}}
\newcommand{\beqar}{\begin{eqnarray}}
\newcommand{\eeqar}{\end{eqnarray}}
\newtheorem{theorem}{Theorem}[section]

\newtheorem{definition}[theorem]{Definition}
\newtheorem{remark}[theorem]{Remark}

\newtheorem{exit}[theorem]{Example}

\newenvironment{rem}{\begin{remark}\rm}{\end{remark}}
\newenvironment{ex}{\begin{exit}\rm}{\end{exit}}

\newcommand{\RR}{{\mathbb R }}
\nc{\FF}{ {\mathbb F} }
\nc{\HH}{ {\mathbb H} }
\newcommand{\ZZ}{{\mathbb Z }}
\newcommand{\QQ}{{\mathbb Q }}
\newcommand{\UU}{{\mathbb U }}

\newcommand{\calf}{\mathcal{F}}

\newcommand{\calo}{\mathcal{O}}
\newcommand{\calp}{\mathcal{p}}

\newcommand{\caly}{\mathcal{Y}}

\newcommand{\senv}{\!\! \smalliint \! \!}

\nc{\umax}{{U_{\max}}}

\newcommand{\liek}{{\mathfrak k}}
\newcommand{\lieu}{{\mathfrak u}}

\newcommand{\hH}{\hat{H}}
\newcommand{\hU}{\hat{U}}

\newcommand{\lieks}{{\liek}^*}
\newcommand{\liet}{{\mathfrak t}}

\newcommand{\xu}{X\env U}

\newcommand{\tkimp}{(T^*K)_{{\rm impl}}}

\nc{\lieq}{{\mathfrak q}}
\nc{\liez}{{\mathfrak z}}
\nc{\lieqs}{{\lieq}^*}
\nc{\lieg}{{\mathfrak g}}
\nc{\liegs}{{\lieg}^*}
\nc{\liep}{{\mathfrak p}}
\nc{\lieps}{{\liep}^*}
\newcommand{\ximpq}{X_{{\rm KimplK^{(P)}}}}
\newcommand{\tkimpq}{(T^*K)_{{\rm KimplK^{(P)}}}}
\newcommand{\tcone}{\liets_{(P)+}}

\def\a{\alpha}
\def\b{\beta}

\def\e{\varepsilon}

\def\l{\lambda}

\def\s{\sigma}

\title{Symplectic quotients of unstable Morse strata for normsquares of moment maps}

\author{Frances Kirwan}

\address{Mathematical Institute\\ Oxford University\\ Woodstock Road, Oxford OX2 6GG\\  UK} 
\email{kirwan@maths.ox.ac.uk}

\begin{document}

\begin{abstract}
Let $K$ be a compact Lie group and fix an invariant inner product on its Lie algebra $\liek$. Given a Hamiltonian action of $K$ on a compact symplectic manifold $X$ with moment map $\mu: X \to \liek^*$,  the normsquare $|\!|\mu |\!|^2$ of $\mu$ defines a Morse stratification $\{ S_\beta: \beta \in \mathcal{B} \}$ of $X$ by locally closed symplectic submanifolds of $X$ such that the stratum to which any $x \in X$ belongs is determined by the limiting behaviour of its downwards trajectory under the gradient flow of $|\!|\mu |\!|^2$ with respect to a suitably compatible Riemannian metric on $X$. The open stratum $S_0$ retracts $K$-equivariantly via this gradient flow to the minimum $\mu^{-1}(0)$ of $|\!|\mu |\!|^2$ (if this is not empty). 
If $\beta \neq 0$  the usual \lq symplectic quotient' $(S_\beta \cap \mu^{-1} (0))/K$ for the action of $K$ on the stratum $S_\beta$ is empty. Nonetheless, motivated by recent results in non-reductive geometric invariant theory, we find that the symplectic quotient construction can be modified to provide natural \lq symplectic quotients' for the unstable strata with $\beta \neq 0$. There is an analogous infinite-dimensional picture for the Yang--Mills functional over a Riemann surface with strata determined by Harder--Narasimhan type.

\end{abstract}

\maketitle

\section{Introduction}

Let $(X,\omega)$ be a compact symplectic manifold, and let $\mu: X \to \lieks$ be a moment map for a Hamiltonian action of a compact Lie group $K$ 
 on $(X, \omega)$. Then the symplectic quotient (or Marsden--Weinstein reduction of $X$ at 0 \cite{MW}) is given by $X \senv K = \mu^{-1}(0)/K$ with its induced symplectic structure. 
Let us fix a $K$-invariant inner product
on the Lie algebra $\liek$ of $K$; 
 then the associated normsquare $|\!|\mu|\!|^2$ of $\mu$ can be considered 
as a Morse function on $X$. This is not in general a Morse function in
the classical sense, nor even a Morse-Bott function, since
the connected components of its set of critical points are 
not in general submanifolds of $X$. Nonetheless, given a suitably compatible $K$-invariant Riemannian metric on $X$, there is
a Morse stratification $\{ S_{\beta} : \beta \in  \mathcal{B}\}$
of $X$ induced by 
$|\!|\mu|\!|^2$ such that each stratum $S_{\beta}$ is a $K$-invariant
locally closed symplectic submanifold of $X$ \cite{K}. Here the stratum to which $x \in X$ belongs is determined by
the limit set of its path of steepest descent for $|\!|\mu|\!|^2$, and the index $\beta$ is the intersection with a positive Weyl chamber $\liet_+$ for $K$ of the co-adjoint orbit which is the image under $\mu$ of the corresponding critical set.


We can attempt to construct symplectic quotients for the restrictions of the Hamiltonian $K$-action to the strata $S_\beta$. However the usual construction, given by $(S_\beta \cap \mu^{-1}(0))/K$, is empty if $\beta \neq 0$. When $K=T$ is abelian we can deal with this problem by shifting the moment map by a suitable constant; a natural choice is to replace $\mu^{-1}(0)$ here with $\mu^{-1}( (1 + \epsilon)\beta) $ for $0 < \epsilon <\!< 1$.  However when $\beta$ is not central $\mu^{-1}( (1 + \epsilon)\beta) $ will not in general be $K$-invariant, so we must modify the construction.
We will see in this paper that this can be done by recalling that there are natural identifications
$$S_\beta \cong K \times_{K\beta} (Y_\beta \cap S_\beta)$$
where $Y_\beta$ is the locally closed submanifold  given for the Morse--Bott function $\mu_\beta(x) = \mu(x).\beta$ by
$$Y_\beta = \{ y \in X : \mbox{the downward trajectory of $y$ for $\text{grad}(\mu_\beta)$ has a limit point $x$ with } \mu_\beta(x) = |\!|\beta|\!|^2 \} , 
$$
and the stabiliser $K_\beta$ of $\beta$ under the adjoint action of $K$ acts diagonally on the product of $K$ with the open subset $Y_\beta \cap S_\beta$ of $Y_\beta$. For sufficiently small $\epsilon > 0$ we will see that  $Y_\beta \cap \mu^{-1}((1 + \epsilon)\beta) \subseteq S_\beta$ is compact and 
$$ S_\beta \, \senv  K := (Y_\beta   \cap \mu^{-1}((1 + \epsilon)\beta)/K_\beta $$
has an induced symplectic structure with which it can be regarded as a symplectic quotient for the $K$-action on the stratum $S_\beta$.

The motivation for this construction comes from the relationship between symplectic quotients and geometric invariant theory (GIT).  
Suppose  that $X \subseteq \PP^n$ is a nonsingular complex projective variety, that $\omega$ is the restriction to $X$ of the Fubini-Study K\"{a}hler form on the complex projective space $\PP^n$ and that $K$ acts linearly on $X$ via a unitary representation $\rho: K \to U(n+1)$. Then the open stratum $S_0$ coincides with the semistable locus $X^{ss}$ in the sense of Mumford's GIT for the induced linear action on $X$ of the complexification $G=K_\CC$ of $K$, and the inclusion $\mu^{-1}(0) \to X^{ss}$ composed with the quotient map from $X^{ss}$ to the GIT quotient $X/\!/G$ induces an identification of the symplectic quotient $X \senv K = \mu^{-1}(0)/K$ with $X/\!/G$. (For this reason, even in the non-algebraic case, we will refer to the strata $S_\beta$ for $\beta \neq 0$ as the unstable strata). The unstable strata $S_\beta$ 
 are $G$-invariant locally closed subvarieties of $X$ and have descriptions of the form $$S_\beta = KY_\beta^{ss} = GY_\beta^{ss} \cong G \times_{P_\beta} Y_\beta^{ss} \cong K \times_{K_\beta} Y_\beta^{ss}$$
where $P_\beta$ is a parabolic subgroup of $G$, and $Y_\beta^{ss} = Y_\beta \cap S_\beta$ has an inductive description involving semistability for the action of a Levi subgroup $L_\beta$ of $P_\beta$, after twisting the linearisation by a suitable rational character of $P_\beta$. 
Since such characters do not in general extend to $G$, in order to construct GIT quotients of the unstable strata $S_\beta$ it is natural to consider quotients of the subvarieties $Y_\beta$ by the action of the non-reductive groups $P_\beta$.
Recent results have extended classical GIT to suitable non-reductive linear algebraic group actions on projective varieties \cite{BDHK,BDHK2}, and the modified symplectic quotient construction for unstable strata just described is suggested by these advances, together with links between non-reductive GIT and the symplectic implosion construction of Guillemin, Jeffrey and Sjamaar \cite{GJS,implone}. In the algebraic setting the modified symplectic quotient construction coincides with a non-reductive GIT quotient construction for the $P_\beta$ action on $Y_\beta$ with an appropriately twisted linearisation.

In their fundamental paper \cite{AB} Atiyah and Bott observed that the Yang--Mills functional over a compact Riemann surface $\Sigma$ plays the role of $|\!|\mu |\!|^2$ in an infinite-dimensional analogue of this picture (modulo a constant which depends on the addition of a central constant to the moment map). Here the corresponding analogue of the GIT or symplectic quotient is a moduli space of semistable holomorphic bundles of fixed rank and degree over $\Sigma$, and the stratification $\{ S_\beta : \beta \in \mathcal{B} \}$
is by the Harder--Narasimhan type of a holomorphic bundle. 
The primary motivation for considering the Yang--Mills functional in \cite{AB} (and $|\!|\mu |\!|^2$  in the finite-dimensional setting explored in \cite{K}) as a Morse function  was to study the cohomology (at least the Betti numbers) of the symplectic quotient. This was done by relating the equivariant cohomology of the compact symplectic manifold $X$, or its infinite-dimensional analogue in the Yang--Mills case, to the equivariant cohomology of the strata, and by describing the equivariant cohomology of the unstable strata inductively in terms of semistable strata for symplectic submanifolds of $X$ acted on by compact subgroups of $K$ (or subgroups of its infinite-dimensional analogue, the relevant gauge group). Later work \cite{Witten,JK,JK2} showed how related ideas could be used to study intersection pairings on $X \senv K$ and the ring structure of its cohomology. In a future paper \cite{BK} we will show how to extend these results to symplectic quotients of unstable strata and other non-reductive GIT quotients.

The layout of this paper is as follows. In $\S$2 we will review the Morse stratification for the normsquare of a moment map on a compact symplectic manifold with a compact Hamiltonian action. In $\S$3 and $\S$4 we will summarise the relevant results from non-reductive GIT and symplectic implosion. Finally $\S$5 describes the construction of symplectic quotients of unstable strata for compact Hamiltonian actions on compact symplectic manifolds, with the main results summarised in Theorem \ref{mainresult}, and $\S$6 considers the infinite-dimensional Yang--Mills analogue.

\section{Normsquares of moment maps and their Morse stratifications}

Suppose that a compact 
Lie group $K$ with Lie algebra ${\liek}$ acts smoothly
on a symplectic manifold
$X$ and preserves the symplectic form $\omega$. 
Any $a\in {\liek}$ determines
a vector field $x\mapsto a_x$ on $X$ defined by
the infinitesimal action of $a$.
A moment map for the action of $K$ on $X$ is a smooth $K$-equivariant map
$\mu :X\rightarrow {\liek}^{\ast}$
which satisfies
$$d\mu(x)(\xi).a=\omega_x(\xi,a_x)$$
for all $x\in X$, $\xi\in T_xX$ and $a\in {\liek}$. Equivalently, 
if $\mu_a:X \to {\RR}$ denotes the component 
of $\mu$ along
$a\in {\liek}$ defined for all $x\in X$ by the pairing
$\mu_a(x)=\mu(x).a$
between $\mu(x) \in {\liek}^{\ast}$ and
$a \in {\liek}$, then $\mu_a$ is a Hamiltonian function                       
for the vector field on $X$ induced by
$a$. 

If the stabiliser $K_{\zeta}$ of $\zeta\in {\liek}^{\ast}$ under the adjoint action of $K$ 
acts with only finite stabilisers on $\mu^{-1}(\zeta)$, then $\mu^{-1}(\zeta)$ is
a submanifold of $X$ and the symplectic form $\omega$ induces a
symplectic structure on the orbifold $\mu^{-1}(\zeta)/K_{\zeta}$ which is the Marsden--Weinstein
reduction at $\zeta$ of the action
of $K$ on $X$. The symplectic quotient $X \senv K$ is the Marsden--Weinstein reduction $\mu^{-1}(0)/K$ at 0. 

The reduction
$\mu^{-1}(\zeta)/K_{\zeta}$ also inherits a symplectic structure
when the action of $K_{\zeta}$
on $\mu^{-1}(\zeta)$ has positive-dimensional stabilisers, but in this case it is likely
to have more serious singularities.

\begin{rem} \label{remalgsit}
Let 
$X$ be a nonsingular complex projective variety
embedded in complex projective space $\PP^n$, and let $G=K_\CC$
be a reductive complex Lie group with maximal compact subgroup $K$ acting on $X$ via a
representation $\rho:G\rightarrow GL(n+1;\CC)$. 
By choosing coordinates on $\PP^n$ appropriately we can assume that $\rho$ maps 
 $K$
 into the unitary group $U(n+1)$. Then the Fubini-Study form $\omega$ on $\PP^n$ restricts to
a $K$-invariant K\"{a}hler form on $X$, and there is a moment map
$\mu :X\rightarrow {\liek}^{\ast}$ defined (up to multiplication by a 
constant scalar factor depending on the convention chosen for
the normalisation of the Fubini-Study form) by
\begin{equation} \mu(x).a = \frac{\overline{\hat{x}}^{t}\rho_{\ast}(a)\hat{x}}
{2\pi i|\!|\hat{x}|\!|^2} \label{mmap} \end{equation}
for all $a\in {\liek}$, where $\hat{x}\in {\CC}^{n+1}-\{0\}$ is a representative
vector for $x\in \PP^n$ and the representation $\rho:K \to U(n+1)$ induces
$\rho_*: \liek \to {\lieu}(n+1)$ and dually $\rho^*:{\lieu}(n+1)^* \to \lieks$.

In this situation  the symplectic quotient 
$X \senv K = \mu^{-1}(0)/K$  coincides with 
the GIT quotient $X/\!/G$  in algebraic geometry described in $\S$3 below \cite{K}. Moreover $x \in X$ lies in the semistable locus $X^{ss}$ if and only if the closure of its orbit $Gx$ meets $\mu^{-1}(0)$, and $x$ lies in the stable locus if and only if its orbit $Gx$ meets the open subset $\mu^{-1}(0)_{\text{reg}}$ of $\mu^{-1}(0)$ where $d\mu$  is surjective. The inclusion $\mu^{-1}(0) \to X^{ss}$ composed with the quotient map $X^{ss} \to X/\!/G$ is $K$-invariant and induces a bijection $\mu^{-1}(0)/K \to X/\!/G$ which can be used to identify the symplectic quotient $X\senv K = \mu^{-1}(0)/K$ with the GIT quotient $X/\!/G$.

When $X$ is K\"ahler but not necessarily algebraic then $\mu^{-1}(0)/K$ inherits a K\"ahler structure (at least away from its singularities) by identifying $\mu^{-1}(0)_{\text{reg}}/K$ with the quotient by $G$ of the open subset $G\mu^{-1}(0)_{\text{reg}}$ of $X$.
\end{rem}

Let us fix a maximal torus $T$ of $K$ and an inner product on the Lie algebra
$\liek$ which is invariant under the adjoint action of $K$, and which we will use to identify $\lieks$ with $\liek$. We will assume that the inner product is chosen so that its restriction to the Lie algebra $\liet$ of $T$  takes rational values on the lattice given by the derivatives at the identity of homomorphisms $S^1 \to T$. Then we can consider the associated normsquare $|\!|\mu|\!|^2$ of $\mu$
as a Morse function on $X$; it is not a Morse function in the classical sense, nor even a Morse--Bott function, but it is shown in \cite{K} that $|\!|\mu|\!|^2$ is a \lq minimally degenerate' Morse function. 
 More precisely,  the set of critical points for $f=||\mu||^2$ is a finite disjoint union of closed subsets $\{C_\beta : \beta \in \mathcal B\}$ along each of which $f$ is \emph{minimally degenerate} in the following sense.

\begin{defn}  A locally closed submanifold $\Sigma_\beta$ containing $C_\beta$ with orientable normal bundle in $X$
 is a \emph{minimising submanifold} for $f=|\!|\mu|\!|^2$ if
  
  \begin{enumerate} \item  the restriction of $f$ to $\Sigma_\beta$ achieves its minimum value exactly on $C_\beta$, and
  
   \item the tangent space to $\Sigma_\beta$ at any point $x \in C_\beta$ is maximal among subspaces of $T_x X$ on which the Hessian $H_x(f)$ is non-negative. 
   
    \end{enumerate} 
    If a minimising submanifold $\Sigma_\beta$ exists, then $f$ is called minimally degenerate along $C_\beta$. 
\end{defn}
    
    In \cite{K} it is shown that as a consequence $|\!|\mu|\!|^2$ induces a smooth stratification $\{S_\beta : \beta \in \mathcal B\}$ of  $X$ such that, for a suitable choice of Riemannian metric (which can be taken to be the K\"ahler metric if $(X,\omega)$ is K\"ahler), $x \in X$ lies in the stratum $S_\beta$ if its path of steepest descent for $|\!|\mu|\!|^2$ has a limit point in the critical subset $C_\beta$.   The stratum $S_\beta$ then coincides with $\Sigma_\beta$ near $C_\beta$.

\begin{rem} \label{almostcx}
Here we choose a $K$-invariant Riemannian metric which is compatible with the symplectic structure in the sense that 
 $X$ has a $K$-invariant  almost-complex structure  such that if $\xi \in T_x X$ then $i\xi$ is the dual with respect to  the metric  of the linear form $\zeta \rightarrow \omega_x (\zeta, \xi)$ on $T_x X$. This implies  that
        $$\text{grad} \mu_\beta (x)  = i \beta_x$$ for all $x \in X$, where $\mu_\beta(x) = \mu(x).\beta$. 
\end{rem}

It is  shown in \cite{K} that  in the situation of Remark \ref{remalgsit} the open stratum $S_0$ of this stratification $\{S_\beta : \beta \in \mathcal B\}$ coincides with the semistable locus $X^{{\text{ss}}}$, and that each stratum $S_\beta$ has the form
     $$S_\beta \cong G \times_{P_\beta}  Y_\beta^{{\text{ss}}}$$
     where $Y^{{\text{ss}}}_\beta$ is a locally closed nonsingular subvariety of $X$ and $P_\beta$ is a parabolic subgroup of $G$ (\cite{K} Theorem 6.18).  Moreover, there is a linear action of a Levi subgroup $L_\beta$ of $P_\beta$ on a nonsingular closed subvariety
     $Z_\beta$ of $X$ such that $Y^{{\text{ss}}}_\beta$ retracts equivariantly  onto the subset $Z^{{\text{ss}}}_\beta$
of semistable points for this action. 
It is also shown in \cite{K} that $|\!|\mu|\!|^2$ is equivariantly perfect, in the sense that its equivariant Morse inequalities are in fact equalities, and that this leads to an {inductive procedure} for  calculating  $\dim H^j_G (X^{{\text{ss}}}; \mathbb Q)$, which in good cases give {the Betti numbers of the quotient variety} $X/\!/G$.  

When $X$ is merely a compact symplectic manifold  acted on by a compact group $K$, the function $||\mu||^2$ still induces a smooth stratification of $X$ and is $K$-equivariantly perfect, providing a formula for the Betti numbers of the symplectic quotient $\mu^{-1}(0)/K$ in the good case when 0 is a regular value of $\mu$, which involves the $K$-equivariant cohomology of the critical subsets $C_\beta$. Indeed the same is true when $||\mu||^2$ is replaced with any convex function of $\mu$  (cf. \cite{AB} \S\S8,12).

 The  set $\mathcal B$ indexing the critical subsets $C_\beta$ and the stratification $\{S_\beta : \beta \in \mathcal B\}$ can be identified with a finite set of orbits of the adjoint representation of  $K$ on its  Lie algebra $\frak k$ (which is identified with its dual using the fixed invariant inner product).  Each orbit in $\mathcal B$ is the image under the moment map $\mu : X \to \frak k^\star \cong \frak k$ of the critical subset which it indexes. If a choice is made of a positive Weyl chamber $\frak t_+$ in the Lie algebra of some maximal torus $T$ of $K$, then each adjoint 
 orbit intersects $\frak t_+$ in a  unique point, so  ${\mathcal B}$  can also be identified with a finite set of points
 in $\frak t_+$. In the situation of Remark \ref{remalgsit} a point of $\frak t_+$ lies in ${\mathcal B}$ if it is the closest point to the origin of the convex hull of a nonempty set of the weights of the unitary representation of $K$ which defines its action on $X \subseteq \PP^n$.  The same is true more generally if we interpret weight here as the image under the $T$-moment map of a connected component of the fixed point set $X^T$.

When ${\mathcal B}$  is identified with a finite set of points
 in $\frak t_+$, for $\beta \in {\mathcal B}$ the submanifold  $Z_\beta$ of $X$  is the union of those components of the fixed point set of the subtorus $T_\beta$ of $K$ generated by $\beta$ on which the moment map for $T_\beta$ given by composing $\mu$ with the restriction map from $\lieks$ to $\liets_\beta$ takes the value $\beta$. Then
$$C_\beta = K(Z_\beta \cap \mu^{-1}(\beta)) \cong K \times_{K_\beta} (Z_\beta \cap \mu^{-1}(\beta)) $$
where the subgroup  $K_\beta$ is the stabiliser of $\beta$ under the adjoint action of $K$ on its Lie  algebra, and in the K\"ahler case, the complexification $L_\beta$ of $K_\beta$ is a Levi subgroup of the parabolic subgroup $P_\beta$ of $G=K_\CC$. 

Since the moment map is $K$-equivariant the image of $Z_\beta$ under $\mu$ is contained in the Lie algebra of $K_\beta$, and thus $\mu |_{Z_\beta}$ can be regarded as a moment map for the action of $K_\beta$ on $Z_\beta$. As moment maps are only determined up to the addition of
a central constant,  $\mu |_{Z_\beta} - \beta$ is also a moment map for the
action of $K_\beta$ on $Z_{\beta}$. 

\begin{rem} In the situation of Remark \ref{remalgsit} this change of moment map  corresponds to a modification of
the linearisation of the action of $K_\beta$ on $Z_{\beta}$, and we define
$Z_{\beta}^{ss}$ to be the set of semistable points of $Z_{\beta}$ with respect
to this modified linear action. Equivalently, $Z_{\beta}^{ss}$ is the stratum
labelled by 0 for the Morse stratification of the function $|\!|\mu - \beta|\!|^2$ on
$Z_{\beta}$. Then
$$Y_{\beta}^{ss} = p_{\beta}^{-1}(Z_{\beta}^{ss}) $$
where $Y_\beta$ and $p_{\beta}:Y_{\beta} \to Z_{\beta}$ are 
given by 
$ \label{pb} p_{\beta}(x) = \lim_{t \to \infty} \exp (-it\beta) x$ and 
$$Y_\beta = \{ y \in X \, |  \, p_\beta(y) \in Z_\beta \}.$$
If $B$ is the Borel
subgroup of $G$ associated to the choice of positive Weyl chamber $\liet_+$ and if
$P_{\beta}$ is the parabolic subgroup $B K_\beta$, then $Y_{\beta}$ and
$Y_{\beta}^{ss}$ are $P_{\beta}$-invariant and we have
$ S_{\beta} = K Y_\beta^{ss} \cong K \times_{K_\beta} Y_\beta^{ss} \cong G \times_{P_{\beta}} Y_{\beta}^{ss}. $
Moreover when $X$ is nonsingular $Y_{\beta}$ is a nonsingular subvariety of $X$ and $p_{\beta}:Y_{\beta}
\to Z_{\beta}$ is a locally trivial fibration whose fibre is isomorphic to
$\CC^{m_{\beta}}$ for some $m_{\beta} \geq 0$.

An element $g$ of $G$ lies in the parabolic subgroup $P_{\beta}$ if and only if
$\exp(-it\beta) g \exp(it\beta)$ tends to a limit in $G$ as $t \to \infty$, and this limit defines
a surjection $q_{\beta}: P_{\beta} \to L_\beta$ such that
$$ p_\beta( gy) = q_\beta (g) p_\beta(y)$$ for each $g \in P_\beta$ and $y \in Y_\beta$.
Since $G=KB$ and $B \subseteq P_{\beta}$
we have $G\overline{Y_{\beta}} = K \overline{Y_{\beta}}$, which is compact, and hence
\begin{equation} \label{closure} \overline{S_{\beta}} \subseteq G\overline{Y_{\beta}}
\subseteq S_{\beta} \cup \bigcup_{|\!| \gamma |\!| > |\!| \beta |\!|} S_{\gamma}. \end{equation}

\end{rem}

\section{Non-reductive geometric invariant theory}

\subsection{GIT for reductive groups}

In Mumford's classical {G}eometric {I}nvariant {T}heory we choose a linearisation  of an action of a   reductive group $G$ on a complex  projective variety  $X$; this  
is given by  an ample line bundle $L$ on $X$ and a lift of the action to $L$.
When $L$ is very ample, so that $X$ can be embedded in a projective space $\PP^n$ such that $L$ is the restriction of the hyperplane line bundle $\calo(1)$,   the action is given by a
representation $\rho:G \to GL(n+1)$ and $ {\hat{\calo}}_L(X) =  \bigoplus_{k= 0}^{\infty} H^0(X,L^{\otimes k})$ is $\kk[x_0,\ldots,x_n]/\mathcal{I}_X$ where $\mathcal{I}_X$ is the ideal generated by the homogeneous polynomials which vanish on $X$.
$$\begin{array}{ccccl}
(X,L) & \leadsto 
 & {\hat{\calo}}_L(X)&=&  \bigoplus_{k= 0}^{\infty} H^0(X,L^{\otimes k})\\
| &&&& \\
| & & \bigcup \!| & &\\
\downarrow & & & & \\
X/\!/G & {\reflectbox{\ensuremath{\leadsto}}}   
 & {\hat{\calo}}_L(X)^G & & \mbox{ algebra of invariants. }
\end{array}$$
Since $G$ is reductive, the algebra of $G$-invariants ${\hat{\calo}}_L(X)^G$ is  {finitely generated} as a graded  algebra and so defines a projective variety  
$X/\!/G = \mbox{Proj}({\hat{\calo}}_L(X)^G)$. The inclusion of ${\hat{\calo}}_L(X)^G$ in ${\hat{\calo}}_L(X)$ determines a rational map $X - - \rightarrow X/\!/G$ which fits into a diagram
$$\begin{array}{rcccl}
 & X & - - \rightarrow & X/\!/G & \mbox{ projective variety}\\
 & \bigcup  & & || & \\
\mbox{semistable} & X^{ss} & \stackrel{\mbox{onto}}{\longrightarrow} & X/\!/G & \\
 & \bigcup  & & \bigcup & \mbox{open}\\
\mbox{stable} & X^s & \longrightarrow & X^s/G &
\end{array} $$     
where $X^s$ and $X^{ss}$ are open subvarieties of $X$, the GIT quotient ${X}/\!/G$ is a categorical quotient for the action of $G$ on $X^{ss}$ via the  $G$-invariant surjective morphism $\phi_G: X^{ss} \to X/\!/G$, and  
$$\phi_G(x) = \phi_G( y) \Leftrightarrow \overline{Gx} \cap \overline{Gy} \cap X^{ss} \neq \emptyset.$$
\begin{rmk} \label{rem2}
A complex Lie group $G$ is reductive if and only if it is the complexification $G=K_\CC$ of a maximal compact subgroup $K$, 
 and then $X/\!/G = \mu^{-1}(0)/K$ for a suitable  moment map $\mu$ for the action of $K$ (see Remark \ref{remalgsit} above). Indeed, recall that 
in this situation the semistable locus $X^{ss}$ coincides with the open stratum $S_0$, while  $x \in X$ lies in  $S_0 = X^{ss}$ if and only if the closure of its orbit $Gx$ meets $\mu^{-1}(0)$, and $x$ lies in the stable locus if and only if its orbit $Gx$ meets the open subset $\mu^{-1}(0)_{\text{reg}}$ of $\mu^{-1}(0)$ where $d\mu$  is surjective. Then the inclusion $\mu^{-1}(0) \to X^{ss}$ composed with the quotient map $X^{ss} \to X/\!/G$ induces an identification of the symplectic quotient $X\senv K = \mu^{-1}(0)/K$ with the GIT quotient $X/\!/G$.

When $X$ is K\"ahler but not necessarily algebraic then 
we can define an equivalence relation $\sim$ on the open stratum $S_0$ by $x \sim y$ if and only if 
$\overline{Gx} \cap \overline{Gy} \cap X^{ss} \neq \emptyset$; if $x$ lies in the open subset $G \mu^{-1}(0)_{\text{reg}}$ of $S_0$ then $x \sim y$ if and only if $y \in Gx$. Then the inclusion $\mu^{-1}(0) \to S_0$  induces an identification of the symplectic quotient $X\senv K = \mu^{-1}(0)/K$ with the topological quotient $S_0/\sim$.
Thus $\mu^{-1}(0)/K$ inherits a stratified K\"ahler structure, with the complex structure induced from $S_0$ and the K\"ahler form given by the symplectic form on $\mu^{-1}(0)/K$ \cite{Huck}.
\end{rmk}

The subsets $X^{ss}$ and $X^s$ of $X$ for a linear action of a reductive group $G$ with respect to an ample linearisation $\mathcal{L}$ are characterised by the Hilbert--Mumford criteria 
\cite[Chapter 2]{GIT}, \cite{New}: 

\begin{propn} 
\label{sss} (i) A point $x \in X$ is semistable (respectively
stable) for the action of $G$ on $X$ if and only if for every
$g\in G$ the point $gx$ is semistable (respectively
stable) for the action of a fixed maximal torus $T$ of $G$.

\noindent (ii) A point $x \in X$ with homogeneous coordinates $[x_0:\ldots:x_n]$
in some coordinate system on $\PP^n$
is semistable (respectively stable) for the action of a maximal 
torus $T$ of $G$ acting diagonally on $\PP^n$ with
weights $\a_0, \ldots, \a_n$ if and only if the convex hull
$$\conv \{\a_i :x_i \neq 0\}$$
contains $0$ (respectively contains $0$ in its interior).
\end{propn}

 The projective GIT quotient $X/\!/G$ contains as an open subset the geometric quotient $X^s/G$ of the stable set $X^s$. When $X$ is nonsingular then the singularities of $X^s/G$ are very mild, since the stabilisers of stable points are finite subgroups of $G$. If $X^{ss} \neq X^s \neq \emptyset$ the singularities of $X/\!/G$ are typically more severe, but $X/\!/G$ has a \lq partial desingularisation'  $\tilde{X}/\!/G$ which ( if $X$ is irreducible and $X^s \neq \emptyset$) is also a projective completion of $X^s/G$ and  is itself  a geometric quotient 
$$\tilde{X}/\!/G = \tilde{X}^{ss}/G$$
by $G$ of an open subset $\tilde{X}^{ss} = \tilde{X}^s$ of a $G$-equivariant blow-up $\tilde{X}$ of $X$ \cite{K2}. 
$\tilde{X}^{ss}$ is obtained from ${X}^{ss}$ by successively blowing up along the subvarieties of semistable points stabilised by reductive subgroups of $G$ of maximal dimension and then removing the unstable points in the resulting blow-up.
Thus for irreducible $X$ we have\\
i) when $X^{ss} = X^s \neq \emptyset$ the GIT quotient $X/\!/G = X^s/G$ is a projective variety which is a geometric quotient of the open subvariety $X^s$ of $X$;\\
ii) when  $X^{ss} \neq X^s \neq \emptyset$ then the GIT quotient $X/\!/G$ is a projective completion of the geometric quotient $X^s/G$, and $X^s/G$ has another projective completion $\tilde{X}/\!/G = \tilde{X}^{s}/G$ which is a \lq partial desingularisation' of $X/\!/G$ in the sense just described.

\begin{rem} The GIT quotient $X/\!/G$ has an ample line bundle which pulls back to a positive tensor power on $X^{ss}$ of the line bundle $L$ defining the linearisation $\mathcal{L}$. 

Note that when we replace the linearisation $\mathcal{L}$ for the action of $G$ on $X$ by any positive tensor power of itself, the stable and semistable loci $X^s$ and $X^{ss}$ and the GIT quotient $X/\!/G$ are unchanged. From a symplectic viewpoint the symplectic form and moment map (and the induced symplectic form on the symplectic quotient) are multiplied by a positive integer, but $\mu^{-1}(0)$ is unchanged. In particular this means that it makes sense to multiply a linearisation by a rational character $\chi/m$ of $G$, where $\chi:G \to \CC^*$ is a character and $m$ is a positive integer: from a GIT perspective we can interpret the result as multiplying the induced linearisation on $L^{\otimes m}$ by the character $\chi$. From a symplectic viewpoint we are adding a central constant to the moment map.
\end{rem}

\begin{ex}
\label{subsec:unstable strata}

As we have seen, associated to the $G$-action on $X$ with linearisation $\mathcal{L}$ and an invariant inner product on the Lie algebra of a maximal compact subgroup $K$ of $G$, there is a stratification (the Morse stratification for $|\!|\mu |\!|^2$)
$$ X = \bigsqcup_{\beta \in \mathcal{B}} S_\beta$$ of $X$ by locally closed subvarieties $S_\beta$, 
indexed by a partially ordered finite subset $\mathcal{B}$ of a positive Weyl chamber for the reductive group $G$,  such that 

 (i) $S_0 = X^{ss}$, 

\noindent and for each $\beta \in \mathcal{B}$

 (ii) the closure of $S_\beta$ is contained in $\bigcup_{\gamma \geqslant \beta} S_\gamma$, and

 (iii) $S_\beta = KY_\beta^{ss} = GY_\beta^{ss} \cong G \times_{P_\beta} Y_\beta^{ss} \cong K \times_{K_\beta} Y_\beta^{ss}$

\noindent where
$P_\beta$ is a parabolic subgroup of $G$ acting on  a projective subvariety $\overline{Y}_\beta$ of $X$ with an open subset $Y_\beta^{ss}$ which is determined by the action of the Levi subgroup $L_\beta$ of $P_\beta$ with respect to a suitably twisted linearisation \cite{Hess,K}.
Here the linearisation $\mathcal{L}$ is restricted to the action of the parabolic subgroup $P_\beta$ over $\overline{Y}_\beta$, and then twisted by the rational character $\beta$ of $P_\beta$.

To construct a quotient by $G$ of (an open subset of) an unstable stratum $S_\beta$, we can study the linear action on $\overline{Y}_\beta$ of the parabolic subgroup $P_\beta$. In order to have a non-empty quotient we must modify the linearisation $\mathcal{L}$, and it is natural to do this by twisting it by a rational character; such a character may not extend to a character of $G$, which is why it makes sense to consider the action on $\overline{Y_\beta}$ of the non-reductive group $P_\beta$. Twisting by $\beta$ (or subtracting $\beta$ from the moment map for the maximal compact subgroup $K_\beta$ of $P_\beta$) gives a categorical quotient $Z_\beta/\!/L_\beta \cong (Z_\beta \cap \mu^{-1}(\beta))/K_\beta \cong C_\beta/K$ for the action of $P_\beta$ on $Y_\beta^{ss}$, or equivalently for the action of $G$ on $S_\beta$ (cf. Remark \ref{pb}), but in general this is far from being a geometric quotient. To have hope of a non-empty open subset of $S_\beta$ with a geometric quotient (when $\beta \neq 0$) one can instead try twisting the action of $P_\beta$ on $\overline{Y_\beta}$ by $(1 + \epsilon)\beta$ where $0 < \epsilon <\!< 1$, or by another perturbation of $\beta$ whose restriction to $T_\beta$ is of this form. 
\end{ex}

\subsection{GIT for non-reductive groups}

Motivated by Example \ref{subsec:unstable strata}, let us consider a complex projective variety $X$ acted on linearly (with ample linearisation $\mathcal{L}$) by a  linear algebraic group $H$ which is not necessarily reductive. Then $H = U \rtimes R$ is the semi-direct product of its unipotent radical $U$ by a reductive subgroup $R$; here $R$ is a Levi subgroup of $H$ and is unique up to conjugation by $H$.

An immediate difficulty arises when trying to extend classical GIT to non-reductive linear algebraic groups $H$; this  is that in general we cannot define a projective {variety} $X/\!/H = \mbox{Proj}({\hat{\calo}}_L(X)^H)$ 
because ${\hat{\calo}}_L(X)^H$ is not necessarily
finitely generated as a graded algebra. However in \cite{BDHK,DK} it is shown that given an $H$-action on $X$ with linearisation $\mathcal{L}$ as above, $X$ has open subvarieties $X^s$
(\lq stable points') and $X^{ss}$ (\lq semistable points') with a geometric quotient $X^s \to X^s/H$ and an \lq enveloping quotient' $X^{ss} \to X\env H$, with 
 a diagram
$$\begin{array}{rcccl}
 & X & & & \\ 
 & \bigcup  & &   & \\   
\mbox{semistable} & X^{ss} & \longrightarrow & X\env H & \\
 & \bigcup  & & \bigcup & \mbox{open}\\
\mbox{stable} & X^s & \longrightarrow & X^s/H &
\end{array} $$
where  {\em if} ${\hat{\calo}}_L(X)^H$ is finitely generated then $X\env H = \mbox{Proj}({\hat{\calo}}_L(X)^H)$ as in the reductive case. 
However  $X\env H$ is not always a projective variety; moreover (even when ${\hat{\calo}}_L(X)^H$ is finitely generated and so $X\env H = \mbox{Proj}({\hat{\calo}}_L(X)^H)$ is a projective variety) the $H$-invariant morphism 
$X^{ss} \to X\env H $ is {not necessarily a categorical quotient}, and its image is not in general a subvariety of $X \env H$, only a constructible subset. A final problem is that  there are  in general no obvious analogues in this situation of the Hilbert--Mumford criteria for (semi)stability.

However non-reductive GIT is better behaved when the unipotent radical $U$ of  $H = U \rtimes R$ is \lq internally graded' in the sense that its Levi subgroup $R \cong H/U$ has a central one-parameter subgroup $\lambda: \CC^* \to R$ whose adjoint action on the Lie algebra of $U$ has only strictly positive weights. 
It is shown in \cite{BDHK2, BDHK3,BK15} that, 
provided that we are willing to twist the 
 linearisation for a linear action of $H$ on a projective variety $X$  by an appropriate (rational) character,  many of the good properties of Mumford's GIT hold. Many non-reductive linear algebraic group actions arising  in algebraic geometry are actions of linear algebraic groups with internally graded unipotent radicals: in particular, any parabolic subgroup of a reductive group has this form, as does the automorphism group of any complete simplicial toric variety \cite{cox}, and  the group of $k$-jets of germs of biholomorphisms of $(\CC^p,0)$ for any positive integers $k$ and $p$ \cite{BK15}.

\begin{ex}
 The automorphism group of the weighted projective plane $\PP(1,1,2)$ 
 with weights 1,1 and 2 is
$\mbox{Aut}(\PP(1,1,2)) \cong R \ltimes U$
where $R \cong  GL(2)$ acting on the two-dimensional weight space with weight 1 is reductive, and 
$U \cong (\CC^+)^3$ is unipotent
 with elements given by $(x,y,z)  \mapsto (x,y,z+\lambda x^2 + \mu xy + \nu y^2)$ 
for $(\lambda,\mu,\nu) \in \CC^3$.
\end{ex}

\begin{definition} \label{defn0.1} Let us call a unipotent linear algebraic group $U$ {\em graded unipotent} 
if there is a homomorphism $\lambda:\CC^* \to Aut(U)$ with the weights of the
 $\CC^*$ action on $Lie(U)$ all {strictly positive}.
For such a homomorphism $\lambda$ let
$$\hat{U} = U \rtimes \CC^* = \{(u,t):u \in U, t \in \CC^*\}$$
be the associated semi-direct product of $U$ and $\CC^*$ with multiplication $(u,t)\cdot(u',t') = (u(\lambda(t)(u')),tt')$. We will say that a linear algebraic group $H=U \rtimes R$ has {\em internally graded unipotent radical} $U$ if the centre $Z(R)$ of $R$ has a one-parameter subgroup $\lambda: \CC^* \to Z(R)$ whose adjoint action grades $U$.

When $L$ is very ample, and so induces an embedding of $X$ in a projective space $\PP^n$, we can choose coordinates on $\PP^n$ such that the action of $\CC^*$ on $X$ is diagonal, given by
$$ t \mapsto \left( \begin{array}{cccc} t^{r_0} & 0 & \ldots & 0\\
 0 & t^{r_1} & \ldots & 0\\
 & & \ldots & \\
 0 & 0 & \ldots & t^{r_n} \end{array} \right) $$
where $r_0 \leq r_1 \leq \cdots \leq r_n$. The {\em lowest bounded chamber} for this linear $\CC^*$-action is the closed interval $[r_0,r_j]$ where  $r_0 = \cdots = r_{j-1} <  r_j \leq \cdots \leq r_n$, with interior the open interval $(r_0,r_j)$, unless the action of $\CC^*$ on $X$ is trivial; when the action is trivial so that  $r_0 = r_1 = \cdots = r_n$ we will say that $[r_0,r_0]$ is the lowest bounded chamber and it is its own interior. Note that in the situation above, if $\CC^*$ acts trivially then so does $U$.
\end{definition}

\begin{thm}[\cite{BDHK2, BDHK3}] \label{firstthm} Let $H=U\rtimes R$ be a linear algebraic group with internally graded unipotent radical $U$ acting linearly on a projective
variety $X$ with linearisation $\mathcal{L}$ on a very ample line bundle $L$. Suppose also that semistability coincides with stability for the unipotent radical $U$, in the sense that 
$$ x \in Z_{{\rm min}}  \Rightarrow {\rm Stab}_U(x) = \{ e \} $$
where $Z_{{\rm min}} $ is the union of those connected components of the fixed point set $X^{\CC^*}$ where $\CC^*$ acts on the fibres of $L^*$ with minimum weight. Then the linearisation for the action of $\hat{U}$ on $X$ can be twisted by a rational character of $\hU$ so that 0 lies in the interior of the lowest bounded chamber for the linear $\CC^*$ action on $X$ and\\ 
(i) the algebras ${\hat{\calo}}_{\mathcal{L}^{\otimes c}}(X)^{\hat{U}} = \oplus_{m=0}^\infty H^0(X,L^{\otimes cm})^{\hat{U}}$ of $\hat{U}$-invariants and 
${\hat{\calo}}_{\mathcal{L}^{\otimes c}}(X)^{H} = \oplus_{m=0}^\infty H^0(X,L^{\otimes cm})^{H}$ of $H$-invariants
 are {finitely generated} for any sufficiently divisible integer $c > 0$,
so that the enveloping quotients $X\env \hat{U} = \mbox{Proj}({\hat{\calo}}_{\mathcal{L}^{\otimes c}}(X)^{\hat{U}})$ and $X\env H = (X\env \hat{U})/\!/(R/\lambda(\CC^*))$
are  projective varieties;\\
(ii) $X^{ss,\hat{U}} = X^{s,\hat{U}} $ and also $X^{ss,H}$ and $X^{s,\hat{U}} $ have Hilbert--Mumford descriptions and $X\env \hat{U} = X^{s,\hat{U}} / \hat{U}$ is a geometric quotient of $X^{s,\hat{U}} $ by $\hat{U}$.
 
Moreover, even when the condition that semistability should coincide with stability for the unipotent radical fails,  
 there is \\
(iii) a projective variety, containing the geometric quotient $X^{s,\hat{U}} /\hat{U}$ as an open subset, which
 is  a geometric quotient  $\tilde{X}^{ss, \hU}/\hat{U}$ by $\hU$ of an open subset $\tilde{X}^{ss, \hU}$ of a $\hat{U}$-equivariant blow-up $\tilde{X}$ of $X$, and \\
(iv) an induced linear action of $R/\lambda(\CC^*)$ on $\tilde{X}^{ss, \hU}/\hat{U}$ whose reductive GIT quotient is a projective variety which contains the geometric quotient $X^{s,H}/H$ as an open subset. 
\end{thm}

\section{Symplectic implosion}

The non-reductive
GIT quotients described in $\S$3 can be studied using symplectic techniques closely
related to the \lq symplectic implosion' construction of Guillemin, Jeffrey and
Sjamaar \cite{GJS,implone}. In this paper  this link will be described for the special case of the unstable strata for the moment map normsquare; in \cite{BK} we will explore more general situations. 

For the original construction \cite{GJS} we suppose that $U$ is a maximal unipotent subgroup of a complex reductive group
$G$ acting linearly (with respect to an ample line bundle $L$) on a complex projective variety 
$X$, and we assume that the linear action of $U$ on $X$ 
extends to a linear action of $G$. Then
the algebra of invariants $\bigoplus_{k \geq 0} H^{0}(X,L^{\otimes k})^U$ is a finitely generated graded algebra and the
 enveloping quotient $\xu$ is the associated projective variety
$\Proj(\bigoplus_{k \geq 0} H^{0}(X,L^{\otimes k})^U)$ \cite{Grosshans2}.
It is shown in \cite{GJS} that  if $K$ is a maximal compact subgroup of $G$, and $X$ is given a suitable $K$-invariant K\"{a}hler
form, then $\xu$
can be identified with the \lq symplectic implosion' or \lq imploded cross-section' $\ximp$ of $X$ by $K$. In this section we will recall this construction and its generalisation  \cite{implone} to the situation when $U$ is  the unipotent radical
of any parabolic subgroup $P$ of  $G$.

As before let $(X,\omega)$ be a symplectic manifold on which a compact connected Lie group $K$
acts with a moment map $\mu:X \to \lieks$ where $\liek$ is the Lie algebra of $K$, and fix an invariant inner product on $\liek$, using it to identify $\lieks$
with $\liek$. Let $T$ be a maximal torus of $K$ with Lie algebra $\liet \subseteq \liek$
and Weyl group $W = N_K(T)/T$, and let $\liets_+ \cong \liets/W \cong \lieks/{\rm Ad}^*(K)$
be a positive Weyl chamber in $\lieks$. The {\em imploded cross-section} \cite{GJS}
of $X$ is then
$$ X_{{\rm impl}} = \mu^{-1}(\liets_+)/\approx $$
where $x \approx y$ if and only if $\mu(x) = \mu(y) = \zeta \in \liets_+$ and
$x = ky$ for some element $k$ of the commutator subgroup $[K_\zeta,K_\zeta]$ of the stabiliser $K_\zeta$  of $\zeta$ under the co-adjoint action of
$K$. If $\Sigma$ is the 
set of faces of $\liets_+$ then $\ximp$ is the disjoint union $$ X_{{\rm impl}}  \ = \  \coprod_{\sigma \in \Sigma} \frac{\mu^{-1}(\sigma)}{[K_\sigma,K_\sigma]}
\ =  \ \mu^{-1}((\liets_+)^\circ) \ \ \ \sqcup
 \coprod_{\begin{array}{c}\sigma \in \Sigma\\  \sigma \neq (\liets_+)^\circ \end{array}
 } \frac{\mu^{-1}(\sigma)}{[K_\sigma,K_\sigma]} $$
where $K_\sigma = K_\zeta$ for any $\zeta \in \sigma$. We give $\ximp$
the quotient topology induced from $\mu^{-1}(\liets_+)$, and it 
inherits a stratified symplectic structure, where the strata are the locally
closed subsets $\mu^{-1}(\sigma)/[K_\sigma,K_\sigma]$. Each such stratum is the symplectic
reduction by the action of $[K_\sigma,K_\sigma]$  of a locally closed symplectic
submanifold
$$X_\sigma = K_\sigma \mu^{-1}( \bigcup_{\tau \in \Sigma, \bar{\tau} \supseteq \sigma} \tau)$$
of $X$; locally near every point $y \in \ximp$ can be identified symplectically with the product of the
stratum containing $y$ and a normal cone \cite{GJS}. The induced action of $T$ on 
$\ximp$ preserves this stratified symplectic structure and has a moment map
$$\mu_{\text{impl}}:\ximp \to \liets_+ \subseteq \liets$$
induced by the restriction of $\mu$ to $\mu^{-1}(\liets_+)$. 
 If $\zeta \in \liets_+$ the symplectic reduction of $\ximp$ at
 $\zeta$ for the action of $T$ is the symplectic reduction of $X$ at $\zeta$ for the action
 of $K$:
$$ \frac{\mu_{\text{impl}}^{-1}(\zeta)}{T} = \frac{\mu^{-1}(\zeta)}{T.[K_\zeta,K_\zeta]} = \frac{\mu^{-1}(\zeta)}{K_\zeta}.
$$
The universal imploded cross-section (or universal symplectic implosion) is the imploded cross-section$$ (T^*K)_{{\rm impl}} = K \times \liets_+ / \approx $$
of the cotangent bundle $T^*K \cong K \times \lieks$ with respect to the $K$-action given by the 
right action of $K$ on itself, with an induced action of $K \times T$ from the left action
of $K$ on itself and the right action of $T$ on $K$. Any other implosion 
$\ximp$ can be constructed as the symplectic quotient of the product $X \times \tkimp$
by the diagonal action of $K$ \cite{GJS}.

The universal symplectic implosion $\tkimp$ is always a complex affine variety and its symplectic structure is
given by a K\"{a}hler form. As in \cite{GJS} we can assume for simplicity that $K$ is semisimple and simply connected; for general compact connected
$K$ one can reduce to this case by considering the product $\tilde{K}$ of the centre of $K$ and the universal cover of its commutator subgroup $[K,K]$, and expressing $K$ as $\tilde{K}/\Upsilon$, where $\Upsilon$ is a finite
central subgroup of $\tilde{K}$.  When $B$ is a Borel subgroup of the complexification $G=K_c$  of $K$ with $G = KB$ and $K \cap B = T$, and  $\umax \leq B$
is the unipotent radical of $B$ (and hence a maximal unipotent subgroup of $G$), then
$\umax$ is a Grosshans subgroup of $G$ \cite{Grosshans}. This means that the quasi-affine variety $G/\umax$
can be embedded as an open subset of an affine variety in such a way that its
complement has complex codimension at least two, and so the algebra 
of invariants $\calo(G)^\umax$ is 
finitely generated. By
\cite{GJS} Proposition 6.8 there is a natural $K \times T$-equivariant identification 
$$\tkimp \cong {\rm Spec}(\calo(G)^\umax)$$
of the canonical affine completion ${\rm Spec}(\calo(G)^\umax)$
of $G/\umax$ with $\tkimp$. It follows that if $X$ is a complex projective
variety on which $G$ acts linearly with respect to an ample
line bundle $L$, and $\omega$ is an associated $K$-invariant K\"{a}hler form on
$X$, then the symplectic quotient $\ximp$ of $X \times \tkimp$ by $K$
can be identified with the non-reductive GIT quotient
$$ X/\!/\umax = {\rm Proj}(\hat{\calo}_L(X)^\umax) \cong (X \times {\rm Spec} (\calo(G)^\umax )))/\!/G
\cong \ximp. $$

Suppose now that $U$ is the unipotent radical of a parabolic subgroup $P$ of the complex reductive
group $G$ with Lie algebra $\liep$. By replacing $P$ with a suitable conjugate in $G$,
we can assume that $P$ contains the Borel subgroup $B$ of $G$ and $U \leq \umax$. Then $P = 
U {L^{(P)}} \cong U \rtimes {L^{(P)}}$, where
the Levi subgroup ${L^{(P)}}$ of $P$ contains the complex maximal torus $T_c$ of $G$,
and we can assume in addition that ${L^{(P)}}$ is the complexification of its intersection
$${K^{(P)}} = {L^{(P)}} \cap K = P \cap K$$
with $K$. There is a  subset $S_P$ of the set $S$ of simple roots such that  $P$ is the 
unique parabolic subgroup of $G$ containing $B$ with the property  that if $\alpha \in S$ then the root space $\lieg_{-\alpha} \subseteq \liep$ if and only if $\alpha \in S_P$.
The Lie algebra of $L^{(P)}$ is generated by the root spaces $\lieg_\alpha$ and $\lieg_{-\alpha}$
for $\alpha \in S_P$ together with the Lie algebra $\liet_c = \liet \otimes_\RR \CC$ of the complexification
$T_c$ of $T$, and  the Lie algebra of $U$ is
$$ \lieu = \bigoplus_{\alpha \in R^+: \lieg_\alpha \not\subseteq {\rm Lie}(L^{(P)})} \lieg_\a $$
where $R^+$ is the set of positive roots for $G$. The Lie algebra of $P$ is
$$ \liep = \liet_c \oplus \bigoplus_{\a \in R(S_P)} \lieg_\a $$
where $R(S_P)$ is the union of $R^+$ with the set of all roots which
can be written as sums of negatives of the simple roots in $S_P$.  
We can decompose $\liek^{(P)} = {\rm Lie} K^{(P)}$ and 
$\liet$ as
$$\liek^{(P)} = [\liek^{(P)},\liek^{(P)}] \oplus \liez^{(P)} \ \ \mbox{ and } \ \ \liet = \liet^{(P)} \oplus
\liez^{(P)}$$
where $[\liek^{(P)},\liek^{(P)}]$ is the Lie algebra of the semisimple part $Q^{(P)} = [K^{(P)},K^{(P)}]$ of $K^{(P)}$, while $\liet^{(P)}$ is the Lie
algebra of the maximal torus $T^{(P)} = T \cap [K^{(P)},K^{(P)}]$ of $Q^{(P)}$, and $\liez^{(P)}$ is
the Lie algebra of the centre $Z(K^{(P)})$ of $K^{(P)}$.

When $U = \umax$ the Iwasawa decomposition
$G = K \ \exp(i\liet) \ \umax$
allows us to identify $G/\umax$ with $K\exp(i\liet)$. More generally
there is a decomposition
\begin{equation} \label{Iwas} G = K \times_{K^{(P)}} P = K \times_{K^{(P)}} \ L^{(P)} U = K \times_{K^{(P)}} \ K^{(P)} \exp(i \liek^{(P)}) U
= K \exp(i\liek^{(P)}) U \end{equation}
giving an identification of $G/U$ with $K \exp(i \liek^{(P)})$.

$U$ is a Grosshans subgroup of $G$ \cite{Grosshans2}, and so the algebra of invariants $\calo(G)^U$ is 
finitely generated and $G/U$ has a canonical affine completion
$$ G/U \subseteq \overline{G/U}^{{\rm a}} = {\rm Spec}(\calo(G)^U) $$
where  the complement of the open subset  $G/U$ of the affine variety $\overline{G/U}^{{\rm a}}$
has complex codimension at least two. Therefore if $G$ acts linearly on  a complex projective
variety $X$ with linearisation $\mathcal{L}$
,  then the algebra of invariants
$$\hat{\calo}_{\mathcal{L}}(X)^U \cong (\hat{\calo}_{\mathcal{L}}(X) \otimes \calo(G)^U)^G$$
is finitely generated, and the associated projective variety
$X/\!/U = {\rm Proj}(\hat{\calo}_{\mathcal{L}}(X)^U)$
is isomorphic to the GIT quotient $ (\overline{G/U}^{{\rm a}} \times X)/\!/G$.
It is shown in \cite{implone} that, just as in the case when $U=\umax$, there is a $K$-invariant K\"{a}hler
form on $\overline{G/U}^{{\rm a}}$ which gives  us an identification of $\xu$ with a symplectic quotient of 
$\overline{G/U}^{{\rm a}} \times X$ by $K$, and thus a symplectic description of $\xu$ generalising
the symplectic implosion construction of \cite{GJS}.

To describe this generalised universal symplectic implosion,
let $\Lambda = \ker(\exp |_{\liet})$ be the exponential lattice in $\liet$, and let
$\Lambda^* = {\rm Hom}_{\ZZ}(\Lambda,\ZZ)$ be the weight lattice in $\liets$, so that 
$\Lambda^*_+ = \Lambda^* \cap \liets_+$ is the monoid of dominant weights. For $\lambda \in \Lambda^*_+$
let $V_{\lambda}$ be the irreducible $G$-module with highest weight $\lambda$, and let
$\Pi = \{\varpi_1, \ldots, \varpi_r \}$
be the set of fundamental weights, forming a $\ZZ$-basis of $\Lambda^*$ and a minimal set of generators
for $\Lambda^*_+$. 
Recall that 
there is an isomorphism of
$G \times G$-modules
\begin{equation} \label{thisiso3}
\calo(G) \cong \bigoplus_{\lambda \in \Lambda^*_+} V_{\lambda} \otimes V_\l^* 
\cong \bigoplus_{\lambda \in \Lambda^*_+} V_{\lambda} \otimes V_{\iota \l} \end{equation}
which restricts to  an isomorphism of $G \times T_c$-modules
\begin{equation} \label{thisiso2}
\calo(G)^\umax \cong \bigoplus_{\lambda \in \Lambda^*_+}  V_\l^{(T)} \otimes V_{\lambda}^* 
\cong \bigoplus_{\lambda \in \Lambda^*_+} V_{\lambda}^*. \end{equation}
where $V^{(T)}_\l$ is the irreducible $T_c$-module with highest weight $\l$.
The graded algebra $\calo(G)^\umax$ is generated 
 by its finite-dimensional vector subspace
 $\bigoplus_{\varpi \in \Pi} V_{\varpi}^*$, which 
gives us a closed
$G \times T_c$-equivariant embedding of $\overline{G/U}_{{\rm max}}^{{\rm a}} = {\rm Spec}(\calo(G)^\umax)$
into the affine space $\bigoplus_{\varpi \in \Pi} V_{\varpi}$. It is shown in \cite{GJS} that $\tkimp$ can be identified with the image of this embedding, equipped with the restriction of
 a flat $K$-invariant K\"{a}hler
structure on $\bigoplus_{\varpi \in \Pi} V_{\varpi}$.

To extend this construction to $\overline{G/U}^{{\rm a}}$ when $U$ is the unipotent radical
of a parabolic subgroup $P$ as above, it is observed in \cite{implone} that 
 $\calo(G)^U$ is generated by the smallest (finite-dimensional)
 ${K^{(P)}}$-invariant subspace of  $\calo(G)$
 which contains 
$\bigoplus_{\varpi \in \Pi} V_{\varpi}^* \cong \bigoplus_{\varpi \in \Pi}
 V_{\varpi}^{(T)} \otimes V_{\varpi}^*.$ Here  ${K^{(P)}}$ acts on $\calo(G)$ via left multiplication on $G$.
Let $E^{{(P)}}$ be the dual of this smallest such ${K^{(P)}}$-invariant subspace $(E^{{(P)}})^*$ 
of $\calo(G)$; 
then $(E^{{(P)}})^*$ is fixed pointwise by $U$, and its inclusion in $\calo(G)^U \subseteq \calo(G)$ induces a closed $ {L^{(P)}}
\times G$-equivariant 
embedding of $\overline{G/U}^{{\rm a}} = {\rm Spec}(\calo(G)^U)$ into the affine space $E^{{(P)}}$.
%
Then
$({E^{(P)}})^*$ decomposes under the action of $K \times {K^{(P)}}$ as a direct sum of irreducible
$K \times {K^{(P)}}$-modules
$$({E^{(P)}})^* = \bigoplus_{\varpi \in \Pi} (V_{\varpi}^{({{P}})})^*$$
where $(V_{\varpi}^{({P})})^*$ is the smallest $K \times {K^{(P)}}$-invariant 
 subspace of $\calo(G)$ containing
$V_{\varpi}^*$. 
Moreover
$(V_{\varpi}^{({{P}})})^* \cong  V_{\varpi}^{{K^{(P)}}} \otimes V_{\varpi}^* $
where $V_{\varpi}^{{K^{(P)}}}$ is the irreducible $K^{(P)}$-module with highest
weight $\varpi$, 
so
$${E^{(P)}} = \bigoplus_{\varpi \in \Pi} V_{\varpi}^{({{P}})}
= \bigoplus_{\varpi \in \Pi}  (V_{ \varpi}^{K^{(P)}})^* \otimes
V_{\varpi} .$$
If 
 $v^{(P)}_{\varpi}$ is the vector in $V_{ \varpi}^{({{P}})} \cong
(V_{\varpi}^{K^{(P)}})^* \otimes V_{\varpi}$ representing the inclusion of 
$V_{\varpi}^{K^{{(P)}}}$ in $V_{\varpi}$ then
the embedding  of $G/U \subseteq \overline{G/U}^{{\rm a}}$ in $E^{(P)}$ induced by the inclusion of $(E^{(P)})^*$ in $\calo(G)^U$ 
takes the identity coset $U$ to $\sum_{\varpi \in \Pi} v_{\varpi}^{(P)}$.
Let
$$V_{\varpi}^{{K^{(P)}}} = \bigoplus_{\lambda \in \Lambda_{\varpi}^*} V_{\varpi,\lambda}^{{K^{(P)}}}$$
be the decomposition of $V_{\varpi}^{{K^{(P)}}}$ into weight spaces with weights
$\lambda \in \liets$ under the action of the maximal torus $T$ of $K^{(P)}$. Then 
$V_{\varpi}^{{{(P)}}}$ decomposes as a $K \times T$-module into a sum of 
irreducible $K \times T$-modules
$$V_{\varpi}^{{{(P)}}} \cong \bigoplus_{\lambda} V_{\varpi} \otimes (V_{\varpi,\lambda}^{{K^{(P)}}})^* $$
and  
$v_{\varpi}^{(P)} = \sum_{\lambda} v_{\varpi,\lambda}^{(P)}$
where $v_{\varpi,\lambda}^{(P)} \in V_{\varpi} \otimes (V_{\varpi,\lambda}^{{K^{(P)}}})^*$
represents the inclusion of $V_{\varpi,\lambda}^{{K^{(P)}}}$ in $V_{\varpi}$. In particular
$v_{\varpi,\varpi}^{(P)}$ is a highest  weight vector for the action of
$K \times K^{(P)}$ on $V^{(P)}_{\varpi}$. 
 The embedding  of $G/U \subseteq \overline{G/U}^{{\rm a}}$ in $E^{(P)}$ induced by the inclusion of $(E^{(P)})^*$ in $\calo(G)^U$ 
takes the identity coset to $\sum_{\varpi \in \Pi} v_{\varpi}^{(P)}$. From the 
decomposition $G = K \exp(i\liek^{(P)}) U$ (\ref{Iwas}) and the
compactness of $K$ it follows that the closure $\overline{G/U}^{{\rm a}}$ of the
$G$-orbit of $\sum_{\varpi \in \Pi} v_{\varpi}^{(P)}$ in $E^{(P)}$ is given by the
$K$-sweep
$$\overline{G/U}^{{\rm a}} = K (\overline{\exp (i \liek^{(P)}) \sum_{\varpi \in \Pi} v^{(P)}_{\varpi}})$$
of the closure in $E^{(P)}$ of the $\exp (i \liek^{(P)})$-orbit of $\sum_{\varpi \in \Pi} v^{(P)}_{\varpi}$.
Similarly the closure 
of the $L^{(P)}$-orbit of $\sum_{\varpi \in \Pi} v^{(P)}_{\varpi}$
 is given by
$K^{(P)} (\overline{\exp (i \liek^{(P)}) \sum_{\varpi \in \Pi} v^{(P)}_{\varpi}})$.

There is a unique $K \times {K^{(P)}}$-invariant Hermitian inner product on ${E^{(P)}} = \bigoplus_{\varpi \in \Pi}
V_{\varpi}^{(P)}$
satisfying $|\!| v_{\varpi,\varpi}^{(P)}|\!| = 1$ for each $\varpi \in \Pi$, which is obtained from
$K$-invariant Hermitian inner products on the irreducible $K$-modules $V_{\varpi}$
and their restrictions to 
$K^{(P)}$-invariant Hermitian inner products on the irreducible $K^{(P)}$-modules $V_{\varpi}^{(P)}$.
This gives ${E^{(P)}}$ a flat
K\"{a}hler structure which is $K \times {K^{(P)}}$-invariant. 
If we identify $(V_{ \varpi}^{K^{(P)}})^* \otimes
V_{\varpi}^{K^{(P)}}$ with ${\rm End}(V_{ \varpi}^{K^{(P)}})$ equipped with the
Hermitian structure 
$\langle A,B \rangle = {\rm Trace}(AB^*)$
in the standard way, then $v^{(P)}_{\varpi}$ is identified with the identity map
in ${\rm End}(V_{ \varpi}^{K^{(P)}})$. 

\begin{defn} \label{defncone} Let $\tcone $ be the cone in $\liets$ given by 
$$\tcone  
= \bigcup_{w \in W^{(P)}} {\rm Ad}^*(w) \liets_+$$
where $W^{(P)}$ is the Weyl group of $Q^{(P)}=[K^{(P)},K^{(P)}]$ (which is a subgroup of
the Weyl group $W$ of $K$). 
\end{defn}


It is shown in \cite{implone} that the restriction to the closure $\overline{\exp(i\liet) \sum_{\varpi \in \Pi} v^{(P)}_{\varpi}}$
of the $\exp(i\liet)$-orbit in $E^{(P)}$ of $ \sum_{\varpi \in \Pi} v^{(P)}_{\varpi}$
of  the moment map $\mu^{E^{(P)}}_T$ for the action of $T$ 
on $E^{(P)}$ is a homeomorphism onto the cone
$\tcone $
in $\liets$. Its inverse 
provides a continuous injection
\begin{equation} \label{calffP} \calf^{(P)} : \tcone  \to \overline{G/U}^{{\rm a}} \subseteq E^{(P)} \end{equation}
such that $\mu_T^{E^{(P)}} \circ \calf^{(P)}$ is the identity on $\tcone $.
Moreover $\overline{\exp(i\liet) \sum_{\varpi \in \Pi} v_\varpi^{(P)}}$
is the union of finitely many $\exp(i\liet)$-orbits, each of the form $$
\calf^{(P)}(\s) = \exp(i\liet) \sum_{\varpi \in \Pi, \lambda \in \Lambda^*_\varpi \cap \bar{\s}} v^{(P)}_{\varpi,\lambda}$$
where $\s$ 
is an open face of $\tcone $. Furthermore
the restriction of  the
$K^{(P)}$-moment map $\mu^{E^{(P)}}:E^{(P)} \to (\liek^{(P)})^*$ to the closure
of the $\exp(i\liek^{(P)})$-orbit in $E^{(P)}$ of $ \sum_{\varpi \in \Pi} v^{(P)}_{\varpi}$
is a homeomorphism from $\overline{\exp(i\liek^{(P)}) \sum_{\varpi \in \Pi} v^{(P)}_{\varpi}}$
onto the closed subset
$$\liek^{(P)*}_+ =  {\rm Ad}^*(K^{(P)})\  \tcone $$
of $\liek^{(P)*}$, and  $\overline{\exp(i\liek^{(P)}) \sum_{\varpi \in \Pi} v_\varpi^{(P)}}$
is the union of finitely many $\exp(i\liek^{(P)})$-orbits which
correspond under this homeomorphism to the open faces of $\liek^{(P)*}_{+}$.

The inverse of $\mu^{E^{(P)}}:\overline{\exp(i\liek^{(P)}) \sum_{\varpi \in \Pi} v^{(P)}_{\varpi}} \to \liek^{(P)*}_+$ gives us a continuous
$K^{(P)}$-equivariant map
$$ \calf^{(P)} : \liek_{+}^{(P)*} \to \overline{G/U}^{{\rm a}} \subseteq E^{(P)} $$ extending (\ref{calffP})
such that $\mu_T^{E^{(P)}} \circ \calf^{(P)}$ is the identity on $\liek^{(P)*}_{+}$. This in turn extends to  a continuous
$K \times K^{(P)}$-equivariant surjection
$$\calf^{(P)} : K \times \liek_{+}^{(P)*} \to \overline{G/U}^{{\rm a}} . $$

\begin{defn} \label{dffn}
If $ \zeta \in \liek^{(P)*}_+ = {\rm Ad}^*({K^{(P)}})\tcone = {\rm Ad}^*({K^{(P)}})\liets_{+}$
let $\zeta = {\rm Ad}^*(k)\xi$ with $k \in K^{(P)}$ and $\xi \in \liets_+$, and let $\s_0$ be the open face of 
$\liets_+$ containing $\xi$. Let $\s_0(P)$ be the open face of $\liets_+$ whose closure is
$$\overline{\sigma_0(P)}
= \{ \zeta \in \liets: \zeta \cdot \a = 0 \mbox{ for all }\a \in R_{\s_0} \setminus R^{(P)} \}$$
where $R$ and $R^{(P)}$ are the sets of roots of $K$ and $K^{(P)}$, and
$ R_{\s_0} 
= \{ \alpha \in R: \zeta \cdot \a = 0 \mbox{ for all }\zeta \in \s_0 \},$ 
so that $\s_0(P)$ is an open subset of the open face containing $\s_0$ of the
cone $\tcone$.
Finally let $K_\zeta(P) = k K_\xi k^{-1}$ where $K_\xi(P) = K_{\s_0(P)}$ is the stabiliser
under the adjoint action of $K$ of any element of $\s_0(P)$.
\end{defn}

\begin{rem}
If $\xi$ lies in the interior of $\tcone$ then $K_\zeta(P) = T$ and 
$[K_\zeta(P),K_\zeta(P)]$ is trivial.

\end{rem}

This leads to  the following definition given in \cite{implone} 
of the  ${K^{(P)}}$-imploded cross-section (or generalised symplectic implosion) $\ximpq$.

\begin{defn} \label{impq}
Let $(X,\omega)$ be a symplectic manifold with a Hamiltonian action of $K$
with  moment map $\mu:X \to \lieks$.
Let
$$ \liek^{(P)*}_+ = {\rm Ad}^*({K^{(P)}})\tcone = {\rm Ad}^*({K^{(P)}})\liets_{+}
  = {\rm Ad}^*({Q^{(P)}})\liets_{+}
\subseteq {\liek^{(P)*}}  \label{lieqsplus} $$
be the sweep of $\liets_{+}$ under the co-adjoint action of ${K^{(P)}}$ on $\lieks$, and let $\Sigma^{(P)}$ be the set of open faces of $\liek^{(P)*}_{+}$.
If $\zeta \in \liek^{(P)*}$ let $K_\zeta(P) $ be as in Definition \ref{dffn}.
The {\em ${K^{(P)}}$-imploded cross-section} (or generalised symplectic implosion)
of $X$ is 
$$ \ximpq = \mu^{-1}(\liek^{(P)*}_+)/\approx_{K^{(P)}} $$
where  $x \approx_{K^{(P)}} y$ if and only if 
$\mu(x) = \mu(y) = \zeta \in 
\liek^{(P)*}_+$
 and
$x = \kappa y$ for some $\kappa \in 
[K_\zeta(P),K_\zeta(P)]$.

The {\em universal ${K^{(P)}}$-imploded cross-section} (or universal generalised symplectic implosion for $K^{(P)} \subseteq K$) is the ${K^{(P)}}$-imploded cross-section
$$ \tkimpq = K \times \liek^{(P)*}_+ / \approx_{K^{(P)}} $$
for the cotangent bundle $T^*K \cong K \times \lieks$ with respect to the $K$-action induced from the 
right action of $K$ on itself. 

\end{defn}

 The map 
$\calf^{(P)} : K \times \liek_{+}^{(P)*} \to \overline{G/U}^{{\rm a}}  $ 
 induces a 
$K \times K^{(P)}$-equivariant homeomorphism 
\begin{equation} \label{imp1}
 \tkimpq = K \times \liek^{(P)*}_+ / \approx_{K^{(P)}} 
 \to
\overline{G/U}^{{\rm a}} \subseteq {E^{(P)}}. \end{equation}
 Moreover under this
identification of $ K \times \liek^{(P)*}_+ / \approx_{K^{(P)}} $ with 
$\overline{G/U}^{{\rm a}} \subseteq {E^{(P)}}$, the moment map
for the action of $K \times K^{(P)}$ on $E^{(P)}$ is induced by
the map $
(K  \times \liek^{(P)*}_+)/\approx_{K^{(P)}} \to
\lieks \times \liek^{(P)*}$
given by 
$$(k,\zeta) \mapsto (Ad^*(k)(\zeta),\zeta)).$$

$\overline{G/U}^{\rm a}$ has an induced $K \times K^{(P)}$-invariant K\"{a}hler
structure as a complex subvariety of $E^{(P)}$; it is stratified by its (finitely many) $G$-orbits, and the $K \times K^{(P)}$-invariant K\"{a}hler structure
 on $E^{(P)}$ restricts to a $K \times K^{(P)}$-invariant symplectic structure on 
 each stratum. Under the homeomorphism 
$ \tkimpq  \to
\overline{G/U}^{{\rm a}} $ of (\ref{imp1}) these strata correspond to the
locally closed subsets 
$$ 
\frac{K \times {\rm Ad}^*(K^{(P)})\sigma 
}{\approx^{K^{(P)}}} \cong K^{(P)} \times_{K_\s \cap K^{(P)}} 
\left( \frac{K \times \sigma 
}{\approx^{K^{(P)}}} \right)$$
$$ \cong K^{(P)} \times_{K_\s \cap K^{(P)}} 
\left( \frac{K \times \sigma 
}{
[K_{\s
(P)
},K_{\s
(P)
}]}  \right)$$
of $\tkimpq$ where $\s \in \Sigma 
$ runs over the open faces of $\liets_+$. 

By construction, when $K$ acts on a symplectic manifold $X$ with moment map $\mu:X \to \lieks$, then
the symplectic quotient of $\overline{G/U}^{{\rm a}} \times X = \tkimpq \times X$
by the diagonal action of $K$ can be identified via $\calf^{(P)}$ with $\ximpq$
(and in particular if $X$ is a projective variety with a linear action of 
the complexification $G$ of $K$, then $\ximpq$ can be identified with 
the GIT quotient of $\overline{G/U}^{{\rm a}} \times X$
by the diagonal action of $G$). 
Thus $\ximpq$ inherits a 
stratified $K \times K^{(P)}$-invariant
symplectic structure
$$  \ximpq =  \bigsqcup_{\sigma \in \Sigma} \frac{\mu^{-1}(\sigma)}{\approx^{K^{(P)}}}$$
\begin{equation} =  \mu^{-1}((\liek^{(P)*}_+)^\circ) \sqcup
 \bigsqcup_{\begin{array}{c}\sigma \in \Sigma\\  \sigma \neq (\liets_+)^\circ \end{array}
 } 
 K^{(P)} \times_{K_\s \cap K^{(P)}} 
\left( \frac{\mu^{-1} (\sigma) }{
[K_{\s(P)},K_{\s(P)}]}  \right) \end{equation}
with strata indexed by the set $\Sigma$ of open faces of $\liets_+$,
which are locally
closed symplectic submanifolds of $\ximpq$.
The induced action of ${K^{(P)}}$ on 
$\ximpq$ preserves this symplectic structure and has a moment map
$$\mu_{\ximpq}:\ximpq \to \liek^{(P)*}_+ \subseteq {\liek^{(P)*}}$$
inherited from the restriction of $\mu$ to $\mu^{-1}(\liek^{(P)}_+)$. 

\begin{rem}  If ${K^{(P)}}=T$ and $\zeta \in \liek^{(P)*}_+$ then 
$K_\zeta(P)=K_\zeta$, and
so $X_{{\rm KimplT}}$ is the standard
imploded cross-section $\ximp$ of \cite{GJS}. On the other hand if ${K^{(P)}}=K$ then 
$K_\zeta(P)$ is conjugate to $T$ and $[K_\zeta(P),K_\zeta(P)]$ is trivial
for all $\zeta \in \liek^{(P)*}_+$,
so  $X_{{\rm KimplK}} = T^*K$.
\end{rem}

\begin{rem} \label{helpfulrem}
When $K$ acts holomorphically on a K\"ahler manifold $X$ with moment map $\mu:X \to \lieks$, then the action of $K$ extends to a holomorphic action of its complexification $G=K_\CC$. Since the generalised symplectic implosion $\ximpq$ is the symplectic quotient of $\overline{G/U}^{{\rm a}} \times X = \tkimpq \times X$
by the diagonal action of $K$, it has an induced K\"ahler structure (cf. Remarks \ref{remalgsit}, \ref{rem2}). The open subset $ \mu^{-1}((\liek^{(P)*}_+)^\circ)$ of $\ximpq$ corresponds to the open subset $(G/U) \times X$
of $\overline{G/U}^{{\rm a}} \times X$, and provides a $K^{(P)}$-invariant slice for the action of $U$ on the open subset $ U\mu^{-1}((\liek^{(P)*}_+)^\circ)$ of $X$. 
Thus if $Y$ is a $P$-invariant complex submanifold of $X$ which meets $\mu^{-1}((\liek^{(P)*}_+)^\circ)$, then $Y \cap \mu^{-1}((\liek^{(P)*}_+)^\circ)$ is a  $K^{(P)}$-invariant slice for the action of $U$ on the open subset $ U(Y \cap \mu^{-1}((\liek^{(P)*}_+)^\circ))$ of $Y$, and its closure in $\ximpq$, which is the image of   $ Y \cap \mu^{-1}(\liek^{(P)*}_+)$ in $\ximpq$, has an induced K\"ahler structure. However the singularities of this closure on the image of the boundary of $ Y \cap \mu^{-1}(\liek^{(P)*}_+)$ are likely to be more serious and harder to describe than in the case when $Y$ is $G$-invariant (or equivalently $K$-invariant).

In the general case when $K$ acts on a symplectic manifold $(X,\omega)$ with moment map $\mu:X \to \lieks$, then we can choose a $K$-invariant almost complex structure which is compatible with $\omega$ as at  Remark \ref{almostcx}. If $Y$ is a $K^{(P)}$-invariant almost complex submanifold of $X$ which is invariant under the induced infinitesimal action of $U$, then just as in the case above the image of   $ Y \cap \mu^{-1}(\liek^{(P)*}_+)$ in $\ximpq$ has an induced $K^{(P)}$-invariant symplectic structure, and almost complex structure, such that it can be regarded as an almost-K\"ahler quotient of $Y$ by the infinitesimal action of $U$.
There is an induced Hamiltonian action of 
$K^{(P)}$ (or any subgroup of $K^{(P)}$)  with moment map $\mu_{\text{Yimpl}}$
induced by the restriction of $\mu$ to $Y \cap \mu^{-1}(\liek^{(P)}_+)$, and we can shift this moment map by any constant in the Lie algebra $\liez^{(P)}$ of the centre $Z(K^{(P)})$ of $K^{(P)}$.
It follows from the definition of $\tcone$ that for a generic choice of $\eta$ in $\liez^{(P)}$ we have $K_\eta = K^{(P)}$ and $\eta \in \mu^{-1}((\liek^{(P)*}_+)^\circ)$, so the symplectic quotient by $K^{(P)}$ is given by
$$ \mu_{\text{Yimpl}}^{-1}(\eta)/K^{(P)} =  (Y \cap \mu^{-1}(\eta))/K^{(P)}.$$
When $Y$ is $K$-invariant this simply recovers for us the symplectic reduction of $Y$ at $\eta$ by the action of $K$, but the viewpoint from symplectic implosion allows us to extend this construction to include submanifolds $Y$ which are not $K$-invariant (cf. \cite{GS}).
\end{rem}

\section{Symplectic quotients of unstable strata}

As before let $X$ be a compact symplectic manifold with a Hamiltonian action of a compact group $K$ with moment map $\mu:X \to \lieks$, and choose a compatible $K$-invariant almost complex structure and  Riemannian metric as at Remark \ref{almostcx}. Fix an invariant inner product on $\liek$ with associated norm.

Let $\{ S_\beta : \beta \in \mathcal{B}\}$ be the Morse stratification for the function $||\mu||^2$. Recall that 
the  set $\mathcal B$ indexing the critical subsets $C_\beta$ for $||\mu||^2$ and the stratification $\{S_\beta : \beta \in \mathcal B\}$ can be identified with a finite subset of  a positive Weyl chamber $\frak t_+$ in the Lie algebra of a maximal torus $T$ of $K$, where a point of $\frak t_+$ lies in ${\mathcal B}$ if it is the closest point to the origin of the convex hull of a nonempty set of the weights for the Hamiltonian action of $K$, and  we interpret \lq weight' as the image under the $T$-moment map of a connected component of the fixed point set $X^T$. Then for $\beta \in {\mathcal B}$ the submanifold  $Z_\beta$ of $X$  is the union of those components of the fixed point set of the subtorus $T_\beta$ of $K$ generated by $\beta$ on which the moment map for $T_\beta$ given by composing $\mu$ with the restriction map from $\lieks$ to $\liets_\beta$ takes the value $\beta$, and
$$C_\beta = K(Z_\beta \cap \mu^{-1}(\beta)) \cong K \times_{K_\beta} (Z_\beta \cap \mu^{-1}(\beta)) $$
where the subgroup  $K_\beta$ is the stabiliser of $\beta$ under the adjoint action of $K$ on its Lie  algebra. 

Recall that $\mu |_{Z_\beta}$ can be regarded as a moment map for the action of $K_\beta$ on $Z_\beta$, and so can $\mu |_{Z_\beta} - \beta$ since $\beta$ is central in $K_\beta$. We can define
$Z_{\beta}^{ss}$ to be the stratum
labelled by 0 for the Morse stratification of the normsquare $|\!|\mu - \beta|\!|^2$ of the moment map 
$\mu |_{Z_\beta} - \beta$ on
$Z_{\beta}$.
For $x \in X$ we let $ p_{\beta}(x)$ be the limit $\lim_{t \to \infty} \exp (-it\beta) x$ of the downward trajectory from $x$ for the Morse--Bott function $\mu_\beta = \mu.\beta$ on $X$,
 and define
$$Y_{\beta}^{ss} = p_{\beta}^{-1}(Z_{\beta}^{ss}) $$
where $Y_\beta$ with $p_{\beta}:Y_{\beta} \to Z_{\beta}$ is 
given by 
$$Y_\beta = \{ y \in X \, |  \, p_\beta(y) \in Z_\beta \}$$
(cf. Remark \ref{pb}).
Then $Y_{\beta}$ and
$Y_{\beta}^{ss}$ are $K_{\beta}$-invariant almost-complex submanifolds of $X$ and we have
$$ S_{\beta} = K Y_\beta^{ss} \cong K \times_{K_\beta} Y_\beta^{ss} . $$
Moreover  $p_{\beta}:Y_{\beta}
\to Z_{\beta}$ is a locally trivial fibration whose fibre is isomorphic to
$\CC^{m}$ for some $m \geq 0$ depending on $\beta$ (and possibly also on the connected component of $Z_\beta$ over which the fibre lies).

The locally closed almost-complex submanifold $Y_\beta$ of $X$ is invariant under the action of the maximal torus  $T$ of $K$, and hence so is its closure $\overline{Y_\beta}$. Therefore by a result of Atiyah \cite{A} (see also \cite{B,K}) the image of $\overline{Y_\beta}$ under the $T$-moment map $\mu_T$ given by composing $\mu$ with  restriction $\lieks \to \liets$ is a convex polytope in $\liets$; indeed it is the convex hull of the (finitely many) images of the $T$-fixed points in $\overline{Y_\beta}$. Thus $\mu_T(\overline{Y_\beta})$ is contained in the half-space in $\liets$ consisting of those $\eta \in \liets$ satisfying $\eta.\beta \geqslant |\!| \beta|\!|^2$;  since by assumption $C_\beta = K(Z_\beta \cap \mu^{-1}(\beta))$ is non-empty, $\beta$ is in fact the closest point to 0 in $\liets \cong \liet$ of this convex hull, and a point $y \in X$ lies in $Y_\beta$ if and only if $\beta$ is the closest point to 0 of the image under $\mu_T$ of its trajectory under the gradient flow of $\mu_\beta$.

Recall that $g \in G=K_\CC$ lies in the parabolic subgroup $P_{\beta}$ if and only if
$\exp(-it\beta) g \exp(it\beta)$ tends to a limit in $G$ as $t \to \infty$, and this limit defines
a surjective homomorphism $q_{\beta}: P_{\beta} \to L_\beta$ 
whose kernel is the unipotent radical $U_\beta$ of $P_\beta$. The chosen almost-K\"ahler structure on $X$ is $K$-invariant, and so by the definition of a moment map the gradient flow of $\mu_\beta$ is given by the vector field $x \mapsto i\beta_x$ where $x \mapsto \beta_x$ is the infinitesimal action of $\beta \in \liet$ on $X$. Thus $Y_\beta$ is invariant under the infinitesimal action of $P_\beta$ on $X$.

This means that we can apply the symplectic implosion construction associated to the unipotent radical $U_\beta$ of $P_\beta$ to $\overline{Y_\beta}$ as in Remark \ref{helpfulrem}, and take a symplectic quotient of the result by the induced Hamiltonian action of the maximal compact subgroup $K^{(P_\beta)} = K_\beta$ of $P_\beta$. As discussed in Remark \ref{helpfulrem}, we can shift the moment map for this induced Hamiltonian action by any constant in the Lie algebra $\liez^{(P_\beta)}$ of the centre $Z(K^{(P_\beta)})$ of $K^{(P_\beta)}$,
and for a generic choice of $\eta$ in $\liez^{(P_\beta)}$ we have $K_\eta = K^{(P_\beta)}$ and  the symplectic quotient by $K^{(P_\beta)} = K_\beta$ is given by
$ (\overline{Y_\beta} \cap \mu^{-1}(\eta))/K_\beta.$

By definition if $\eta = \beta$ or any nonzero scalar multiple of $\beta$ then $K_\eta = K_\beta=K^{(P_\beta)} $ and this symplectic quotient is $ (\overline{Y_\beta} \cap \mu^{-1}(\eta))/K_\beta.$ It follows from the description of $Y_\beta$ above that if $\eta$ is a generic element of $\liez^{(P_\beta)}$ and $\eta.\beta < |\!| \beta|\!|^2$ (or if we have equality and $\eta \neq \beta$) then $ (\overline{Y_\beta} \cap \mu^{-1}(\eta))/K_\beta$ is empty, while if $\eta = \beta$ then  the symplectic quotient $ (\overline{Y_\beta} \cap \mu^{-1}(\eta))/K_\beta =   ({Z_\beta} \cap \mu^{-1}(\beta))/K_\beta = C_\beta/K$ collapses the stratum onto its critical subset. On the other hand if $\eta$ is a sufficiently small perturbation of $\beta$ in  $\liez^{(P_\beta)}$ then 
$ \overline{Y_\beta} \cap \mu^{-1}(\eta) \subseteq {Y_\beta}$ so
$$ (\overline{Y_\beta} \cap \mu^{-1}(\eta))/K_\beta = ( {Y_\beta} \cap \mu^{-1}(\eta))/K_\beta.$$
It is therefore natural to choose $\eta$ to be
$(1 + \epsilon)\beta$ for some sufficiently small $\epsilon > 0$ and define 
\begin{equation} \label{defnsqus}
S_\beta \senv_\epsilon K :=  ( {Y_\beta} \cap \mu^{-1}((1 + \epsilon)\beta))/K_\beta.
\end{equation}
This has a stratified symplectic structure and 
it follows from the theory of variation of symplectic quotients \cite{DH,GSt,LS} (cf. \cite{dh98,Thaddeus}) that $S_\beta \senv_\epsilon K$ is independent of $\epsilon$ up to diffeomorphism for $0 < \epsilon <\!< 1$ and the induced symplectic structure varies in a predictable fashion with $\epsilon$; we can also use this theory to study the variation if $\eta$ is chosen to be a different perturbation of $\beta$.

We have thus proved our main result.

\begin{thm} \label{mainresult}
Let $X$ be a compact symplectic manifold with a Hamiltonian action of a compact group $K$ with moment map $\mu:X \to \lieks$. Choose a compatible $K$-invariant almost complex structure and  Riemannian metric on $X$, and an invariant inner product on $\liek$ with associated norm. 
Let $\{ S_\beta : \beta \in \mathcal{B}\}$ be the Morse stratification for the function $||\mu||^2$. If $\beta \in \mathcal{B} \setminus \{0\} $ and $0< \epsilon <\!< 1$ then 
$$
S_\beta \senv_\epsilon K =  ( {Y_\beta} \cap \mu^{-1}((1 + \epsilon)\beta))/K_\beta$$
is a compact stratified symplectic space which can be interpreted as a symplectic quotient for the action of $K$ on the stratum $S_\beta$.

When $X \subseteq \PP_n$ is a complex projective variety equipped with the Fubini--Study K\"ahler form and a linear action of $K$ defined by a unitary representation $K \to U(n+1)$, then when $\epsilon$ is rational $
S_\beta \senv_\epsilon K$ can be identified with a quotient of $S_\beta$ by $G=K_\CC$ obtained from non-reductive GIT applied to the action of the parabolic subgroup $P_\beta$ on $Y_\beta$, with the linearisation twisted by the rational character $(1+\epsilon)\beta$ of $P_\beta$.
\end{thm}


\section {The Yang--Mills functional over a compact Riemann surface}

Atiyah and Bott observed \cite{AB} that the Yang--Mills functional over a compact Riemann surface $\Sigma$ plays the role of $|\!|\mu |\!|^2$ for an infinite-dimensional Hamiltonian action.  Here the symplectic quotient can be identified with a moduli space of semistable holomorphic bundles of fixed rank and degree over $\Sigma$, and the stratification $\{ S_\beta : \beta \in \mathcal{B} \}$
is by the Harder--Narasimhan type of a holomorphic bundle. 

Let $\Sigma$ be a compact Riemann surface of genus $g \geqslant 2$, and let $\mathcal{E}$ be a fixed $C^{\infty}$ complex hermitian vector bundle of rank $n$ and degree $d$ over $\Sigma$. Let $\mathcal{C}$ be the space of all holomorphic structures on $\mathcal{E}$. Since $\Sigma$ has complex dimension one there are no integrality conditions to be satisfied, so $\mathcal{C}$ can be identified with the space of unitary connections on $\mathcal{E}$, which is an infinite-dimensional complex affine space with a flat K\"ahler structure.  

Let $\mathcal{G}_\CC$ denote the group of all $C^{\infty}$ complex automorphisms of $\mathcal{E}$. We can regard $\mathcal{G}_\CC$ as the complexification of the gauge group $\mathcal{G}$ consisting of $C^{\infty}$ unitary automorphisms of $\mathcal{E}$. The natural action of $\mathcal{G}_\CC$ on $\mathcal{C}$ preserves its complex structure and the action of the gauge group $\mathcal{G}$ preserves the K\"ahler structure and is Hamiltonian with a moment map given by the curvature of a connection. The central subgroup $\CC^*$ of $\mathcal{G}_\CC$ given by scalar multiplication on $\mathcal{E}$ acts trivially on $\mathcal{C}$, so the moment map in the direction of the corresponding central $S^1$ in $\mathcal{G}$ is constant; it is essentially given by the ratio $d/n$. The Yang--Mills functional on $\mathcal{C}$ takes a connection to the normsquare of its curvature and hence plays the role of $|\!|\mu|\!|^2$ for the action of the gauge group on $\mathcal{C}$, except that it is more natural to choose the moment map $\mu$ so that $\mu^{-1}(0)$ is nonempty by adding a suitable central constant to the curvature. This means that the Yang--Mills functional differs from $|\!|\mu|\!|^2$ by a constant, so their Morse stratifications will coincide.

Atiyah and Bott \cite{AB} identified the symplectic quotient $\mu^{-1}(0)/\mathcal{G}$ of $\mathcal{C}$ by the gauge group with the 
moduli space $\mathcal{M}(n,d)$ of semistable holomorphic vector bundles of rank $n$ and degree $d$ on $\Sigma$ (modulo S-equivalence).
Recall that a holomorphic vector bundle $E$ over $\Sigma$ is  semistable (respectively stable) if every holomorphic subbundle $D$ of $E$ satisfies
$\text{slope} (D) \leq \text{slope} (E)$, (respectively $\text{slope} D) < \text{slope} (E)$),
where $\text{slope} (D) = \text{deg}(D)/\text{rank}(D)$ 
 (and thus semistable bundles of coprime rank and degree are stable). Any semistable vector bundle $E$ has a Jordan--H\"older filtration by sub-bundles of the same slope as $E$ whose successive subquotients are stable; its associated graded bundle is the direct sum of these successive subquotients (which is independent of the choice of Jordan--H\"older filtration), and two semistable bundles of rank $n$ and degree $d$ are S-equivalent if their associated graded bundles are isomorphic. 

A holomorphic bundle $E$ over $\Sigma$ of rank $n$ and degree $d$ has a canonical Harder--Narasimhan filtration
$$0 = E_0 \subset E_1 \subset \cdots \subset E_{s-1} \subset E_s = E$$
such that $\text{slope} (E_{j-1}) > \text{slope} (E_{j})$ and $E_j/E_{j-1}$ is semistable for $1 \leq j \leq s$. The Harder--Narasimhan type of $E$ is then given by the data provided by the ranks and degrees of the successive subquotients $E_j/E_{j-1}$; in \cite{AB} this is encoded in the vector
$$\lambda(E) = (d_1/n_1, \ldots, d_1/n_1, d_2/n_2, \ldots, d_s/n_s)$$
in which $d_j/n_j$ occurs $n_j$ times.

The main aim of \cite{AB} is to study the cohomology of the moduli space $\mathcal{M}(n,d)$ by showing that the Yang--Mills functional is equivariantly perfect as a Morse function. Because of the analytical difficulties created by working in infinite-dimensions and the singularities in the critical locus for the Yang--Mills functional, Morse theory is not applied directly to the Yang--Mills functional in \cite{AB} but instead the stratification is defined directly in terms of Harder--Narasimhan types; however the analytical difficulties were later overcome \cite{Dask}. Let $\Lambda$ denote the set of all Harder--Narasimhan types, and for any Harder--Narasimhan type $\lambda$, let $\mathcal{C}_\lambda$ denote the subset of $\mathcal{C}$ consisting of holomorphic structures on $\mathcal{E}$ with Harder--Narasimhan type $\lambda$. Atiyah and Bott showed that $\{ \mathcal{C}_\lambda : \lambda \in \Lambda \}$ is a $\mathcal{G}$-equivariantly perfect stratification of $\mathcal{C}$. They conjectured that it coincides with the Morse stratification for the Yang--Mills functional, which was later proved by Daskalopoulos \cite{Dask}.

The moduli space $\mathcal{M}(n,d)$ can also be constructed as  finite-dimensional symplectic or GIT quotients, and the inductive formulas of \cite{AB} for its Betti numbers can be rederived via these \lq finite-dimensional  approximations' to the Yang--Mills picture \cite{ADK,Karxiv}. In \cite{hok12} it is shown that the moduli spaces $\mathcal{M}(n,d)$ (and more generally moduli spaces of sheaves over any fixed nonsingular projective scheme) can be constructed as GIT quotients for actions of complex reductive groups on finite-dimensional complex varieties such that, for any given Harder--Narasimhan type $\lambda$ for bundles of rank $n$ and degree $d$ there is a choice of GIT construction of $\mathcal{M}(n,d)$ for which the bundles of Harder--Narasimhan type $\lambda$ appear as a stratum in the associated stratification $\{S_\beta:\beta \in \mathcal{B}\}$. The results of non-reductive GIT  described in $\S$3 can then be used to construct moduli spaces of holomorphic bundles of fixed Harder--Narasimhan type \cite{behjk16, bhk17}. 

Alternatively we can attempt to use the infinite-dimensional Yang--Mills construction of the moduli space $\mathcal{M}(n,d)$ as a symplectic quotient of $\mathcal{C}$ by the gauge group and the methods of this paper to find an analogous symplectic construction of moduli spaces of holomorphic bundles of fixed Harder--Narasimhan type. Ignoring the analytical difficulties associated to working with infinite-dimensional spaces and groups, we might proceed as follows.

Let $\lambda(E) = (d_1/n_1, \ldots, d_1/n_1, d_2/n_2, \ldots, d_s/n_s)$ be a Harder--Narasimhan type and fix a $C^{\infty}$ filtration 
\begin{equation}
\label{filt}
0 = \mathcal{E}_0 \subset \mathcal{E}_1 \subset \cdots \subset \mathcal{E}_{s-1} \subset \mathcal{E}_s = \mathcal{E}  \end{equation}
of the $C^{\infty}$ bundle $\mathcal{E}$ with $\text{deg}(\mathcal{E}_j/\mathcal{E}_{j-1}) = d_j$ and 
$\text{rank}(\mathcal{E}_j/\mathcal{E}_{j-1}) = n_j$ for $1 \leq j \leq s$. Define $\mathcal{Y}_\lambda$ to be the subset of $\mathcal{C}$ consisting of those holomorphic structures (or equivalently unitary connections) on $\mathcal{E}$ which are compatible with this filtration, in the sense that the subbundles $\mathcal{E}_j$ are all holomorphic subbundles, and define $\mathcal{Y}_\lambda^{ss}$ to consist of those holomorphic structures for which in addition the induced holomorphic structures on the subquotients $\mathcal{E}_j/ \mathcal{E}_{j-1}$ are semistable, so that the holomorphic structure on $\mathcal{E}$ lies in $\mathcal{C}_\lambda$. Let $\mathcal{P}_\lambda$ be the subgroup of $\mathcal{G}_\CC$ consisting of the complex $C^{\infty}$-automorphisms of $\mathcal{E}$ which preserve the filtration (\ref{filt}) and  let $\mathcal{U}_\l$ be the kernel of its induced action on the direct sum of the successive subquotients $\mathcal{E}_j/\mathcal{E}_{j-1}$. There is a  $C^{\infty}$ decomposition of $\mathcal{E}$ as  the orthogonal direct sum of the successive subquotients $\mathcal{E}_j/\mathcal{E}_{j-1}$;  let 
$\mathcal{L}_\lambda$ be the subgroup of $\mathcal{P}_\lambda$ preserving this direct sum decomposition and let $\mathcal{K}_\l$ be its intersection with the gauge group $\mathcal{G}$. Finally let $\mathcal{Z}_\l$ be the subset of $\mathcal{Y}_\l$ consisting of holomorphic structures for which this orthogonal direct sum decomposition of $\mathcal{E}$ is a holomorphic decomposition, and let $\mathcal{Z}_\l^{ss} = \mathcal{Z}_\l \cap \mathcal{Y}_\l^{ss}$. 

Then $\mathcal{Y}_\lambda$, $\mathcal{Y}_\lambda^{ss}$, $\mathcal{Z}_\lambda$, $\mathcal{Z}_\lambda^{ss}$ and $\mathcal{P}_\lambda$, $\mathcal{U}_\lambda$, $\mathcal{L}_\lambda$, $\mathcal{K}_\lambda$ play the roles for the Hamiltonian action of the gauge group $\mathcal{G}$ on $\mathcal{C}$, and on its stratum $\mathcal{C}_\l$, which $Y_\beta$, $Y_\beta^{ss}$, $Z_\beta$, $Z_\beta^{ss}$ and $P_\beta$, $U_\beta$, $L_\beta$ and $K_\beta$ play in the finite-dimensional setting for the Hamiltonian action of the compact group $K$ on the compact symplectic (or K\"ahler) manifold $X$, and its stratum $S_\beta$. Note however that it is really $\lambda - (d/n, \ldots, d/n)$, not $\l$ itself, which plays the role of $\beta$, since the central circle subgroup in the gauge group acts trivially, so
$$(1 + \epsilon) \l - \e (d/n, \ldots, d/n)$$
plays the role of $(1 + \e )\b$.
 
Thus by analogy with the finite-dimensional situation we expect the stratum $\mathcal{C}_\l$ to have a symplectic quotient
$$ \mathcal{C}_\l \senv_\e \mathcal{G} =  ( \mathcal{Y}_\l \cap \text{curv}^{-1}((1 + \epsilon)\l - \e (d/n, \ldots ,d/n)))/\mathcal{K}_\l$$
for $0 < \e <\!< 1$, where $\text{curv}$ assigns to a holomorphic structure, or equivalently a unitary connection, on $\mathcal{E}$ its curvature, appropriately normalised, and $\mathcal{Y}_\l$ and $\mathcal{K}_\l$ are defined as above. At least away from its singularities we expect this symplectic quotient to be identifiable with a suitable moduli space of holomorphic bundles of Harder--Narasimhan type $\l$.

\end{document}